\documentclass[a4paper]{amsart}

\usepackage{amsmath,amssymb,amsfonts,exscale,url,doi,graphicx,xcolor}
\usepackage[ruled,vlined,algosection]{algorithm2e}
\usepackage[numbers,sort]{natbib}

\definecolor{color1}{rgb}{0,0,1}
\hypersetup{pdfborder={1 0 0}}
\renewcommand{\ldots}{\ensuremath{\dotsc}}
\newcommand{\comment}[1]{}

\newcommand{\C}{\mathbb{C}}
\newcommand{\R}{\mathbb{R}}

\newcommand{\X}{\mathcal{X}}
\newcommand{\Y}{\mathcal{Y}}

\newcommand{\K}{\mathcal{K}}
\newcommand{\U}{\,\mathcal{U}}
\newcommand{\V}{\mathcal{V}}
\newcommand{\VP}{{\mathcal{V}_\perp}}
\newcommand{\QVP}{V_\perp}
\newcommand{\Z}{\mathcal{Z}}

\def\la{\lambda}

\def\span{{\rm span}}
\def\dim{{{\rm dim}}}
\def\max{{{\rm max}}}
\def\min{{{\rm min}}}

\newcommand{\diag}{{\rm{diag}}}

\newtheorem{theorem}{Theorem}[section]
\newtheorem{lemma}[theorem]{Lemma}

\newtheorem{remark}{Remark}[section]

\numberwithin{equation}{section}

\begin{document}

\title[Majorization bounds for block filters and eigensolvers]
{Majorization-type cluster robust bounds\\
 for block filters and eigensolvers}

\author[M.~Zhou]{Ming Zhou}
\address{Universit\"at Rostock, Institut f\"ur Mathematik, 
        Ulmenstra{\ss}e 69, 18055 Rostock, Germany}
\email{ming.zhou at uni-rostock (dot) de}

\author[M.~E.~Argentati]{Merico E.~Argentati}
\address{Department of Mathematical and Statistical Sciences,
University of Colorado Denver, Denver, CO 80217-3364, USA}
\email{merico.argentati at ucdenver (dot) edu}

\author[A.~V.~Knyazev]{Andrew V.~Knyazev}
\address{Department of Mathematical and Statistical Sciences,
University of Colorado Denver, Denver, CO 80217-3364, USA}
\email{andrew.knyazev at ucdenver (dot) edu}

\author[K.~Neymeyr]{Klaus Neymeyr}
\address{Universit\"at Rostock, Institut f\"ur Mathematik, 
        Ulmenstra{\ss}e 69, 18055 Rostock, Germany}
\email{klaus.neymeyr at uni-rostock (dot) de}


\keywords{Majorization, principal angles, Ritz values,
 subspace eigensolvers, cluster robustness, block Lanczos method.
 \hfill December 2021}

\begin{abstract}
Convergence analysis of block iterative solvers for Hermitian eigenvalue problems
and the closely related research on properties of matrix-based signal filters
are challenging, and attract increasing attention due to their recent applications
in spectral data clustering and graph-based signal processing.
We combine majorization-based techniques pioneered for investigating
the Rayleigh-Ritz method in [SIAM J. Matrix Anal. Appl., 31 (2010), pp. 1521--1537]
with tools of classical analysis of the block power method by Rutishauser
[Numer. Math., 13 (1969), pp. 4--13] to derive convergence rate bounds
of an abstract block iteration, wherein tuples of tangents of principal angles
or relative errors of Ritz values are bounded using majorization in terms of
arranged partial sums and tuples of convergence factors. 
Our novel bounds are robust in presence of clusters of eigenvalues,
improve some previous results, and are applicable to
most known block iterative solvers and matrix-based filters,
e.g., to block power, Chebyshev, and Lanczos methods
combined with shift-and-invert approaches and polynomial filtering.
\end{abstract}

\maketitle
\pagestyle{myheadings}
\thispagestyle{plain}

\section{Introduction}\label{intro}

Matrix eigenvalue problems have historically appeared in computational mechanics
and then in quantum physics as a tool to approximate solutions of wave and other
time-dependent differential equations. Recent applications in machine learning,
e.g., for spectral clustering and image semantic segmentation, often require
numerical solution of eigenvalue problems for matrices of large sizes and
in some cases in real time, e.g., for autonomous driving. Only a fraction
of the eigenvectors is of practical interest, where the eigenvectors
need to be computed just with certain accuracy. 
Such computations of specific invariant subspaces can be efficiently performed
by iterative eigenvalue problem solvers (eigensolvers)
that simultaneously iterate several vectors in a block.
Recently graph-based signal processing has gained attention.
Therein one utilizes signal filters that are functions, often a polynomial
like Chebyshev, of a matrix with a goal of amplifying the components
of a signal corresponding to selected eigenvalues of the matrix.
Such filters can be interpreted as rudimentary eigensolvers with a fixed accuracy.
Block filters and solvers
match well with modern hardware that can optimize computing matrix-matrix products,
which commonly is the main operation in block algorithms.
 
The convergence speed of block iterative solvers,
or the reduction/amplification quality of signal filters,
determines the efficiency of calculations
and needs to be estimated in advance a priori to bound the computational time
for streaming data processing. Deriving a priori convergence rate bounds
is a classical area of research. The convergence speed of block iterative eigensolvers
is determined by the quality of an initial approximation to target eigenvectors
and the distribution of eigenvalues of the given matrix. Block iterative eigensolvers
are known to be robust with respect to possible clustering of target eigenvalues
-- this property is called ``cluster robustness''. Sharp bounds of convergence rates
of these eigensolvers thus also need to be cluster robust,
which is the main topic of our investigation.
 
Majorization techniques appear naturally for cluster robust bounds involving
eigenvalue approximations \cite{ka10,msz2020}. We explore relations between
Ritz value errors and principal angles in the context of subspace iterates
of block iterative eigensolvers. We refine the majorization-type analysis of tangents
of principal angles from \cite{ka10} by using auxiliary vectors from the classical analysis
of the block power method by Rutishauser \cite{r1969}.
Our new results applied to the block Lanczos method improve existing ones
from \cite{s1980,lz2015}. The improvement is evident
due to the cluster robustness of our majorization based approach.

\subsection{Known cluster robust error bounds}

The cluster robustness is a typical feature of block iterations for solving
matrix eigenvalue problems. A simple example is the block power method
that aims at the largest magnitude eigenvalues of a matrix $A\in\C^{n \times n}$
and significantly improves the simple power method
in the case where the target eigenvalues are clustered.
The convergence of the subspace iterates $\Y^{(0)},\,\Y^{(1)},\,\ldots$
given by $\Y^{(\ell+1)}=A\Y^{(\ell)}$
toward an invariant subspace can be analyzed in terms of various angle measures.
Classical bounds by Rutishauser \cite{r1969} and Parlett \cite[Chapter 14]{p}
can be extended to a normal $A$ with
the eigenvalue arrangement $|\la_1|\ge\cdots\ge|\la_n|$, namely,
\begin{equation}\label{bpme}
 \tan\angle(x_i,\Y^{(\ell)})\le\left|\frac{\la_{p+1}}{\la_i}\right|^{\ell}
 \tan\angle(\X,\Y^{(0)}),\quad i=1,\ldots,p.
\end{equation}
Therein $\angle(\cdot,\cdot)$ denotes Euclidean angles, and the separation
of targeted eigenvalues from the rest is described by the assumption
$|\la_p|>|\la_{p+1}|$ for the block size $p=\dim\Y^{(0)}$. Moreover,
$\X$ denotes the invariant subspace spanned by orthonormal eigenvectors
$x_1,\ldots,x_p$ associated with the first $p$ eigenvalues. In the nontrivial case
$\angle(\X,\Y^{(0)})<\pi/2$, all further subspace iterates have the same dimension $p$
since the first $p$ eigenvalues are nonzero due to $|\la_p|>|\la_{p+1}|\ge0$.
Evidently, the convergence factor $|\la_{p+1}/\la_i|$ is robust with respect to
possible clustering of eigenvalues indexed by $i \le p$.
In contrast, the convergence factor for computing $x_i$ with the simple power method
is $|\la_{i+1}/\la_i|$ which is unimprovable
and might be close to $1$. If $i>1$, a deflation with respect to
the first $i{-}1$ eigenvalues is required.

An essential argument by Rutishauser \cite{r1969} adapted to
\eqref{bpme} is that $\Y^{(\ell)}$ contains
a nonzero vector which is orthogonal to an invariant subspace associated with
$\la_1,\ldots,\la_{i-1},\la_{i+1},\ldots,\la_p$.
Thus $\la_i$ and $\la_{p+1}$ become the only two relevant eigenvalues
in the derivation. This argument is directly applicable to the abstract block iteration
\begin{equation}\label{abi}
 \Y'=f(A)\Y
\end{equation}
that describes various practical methods.
Since $A$ and $f(A)$ share the same eigenvectors,
one attempts to select a proper function $f(\cdot)$ such that
the desired eigenvalues of $A$, if denoted by $\la_1,\ldots,\la_p$, get mapped into
the top $p$ magnitude eigenvalues $f(\la_1),\ldots,f(\la_p)$ of $f(A)$.
Then \eqref{bpme} turns into
\begin{equation}\label{bpmf}
 \tan\angle(x_i,\Y')\le\left|\frac{f(\la_{p+1})}{f(\la_i)}\right|\,
 \tan\angle(\X,\Y),\quad i=1,\ldots,p
\end{equation}
analogously to \cite[Lemma 2.3.1]{k}.
One can use a shifted Chebyshev polynomial $f(\cdot)$
in order to investigate Chebyshev signal filters and the block Lanczos method.
Let eigenvalues of a Hermitian matrix $A$ be arranged as $\la_1\ge\cdots\ge\la_n$
with the assumption $\la_p>\la_{p+1}$, and $f(\cdot)$ be defined by
\begin{equation}\label{fsp}
 f(\alpha)=T_{k-1}\left(1+2\,\frac{\alpha-\la_{p+1}}{\la_{p+1}-\la_n}\right)
\end{equation}
with the Chebyshev polynomial $T_{k-1}$ of the first kind
with degree $k{\,-\,}1$. Then
$|f(\la_i)|>1\ge|f(\la_j)|$ holds for arbitrary indices $i \le p$ and $j>p$
so that the values $f(\la_1),\ldots,f(\la_p)$ are indeed
the top $p$ magnitude eigenvalues of $f(A)$.
Thus, substituting \eqref{fsp}
in \eqref{bpmf} immediately leads to
\begin{equation}\label{blme}
 \tan\angle(x_i,\Y')\le\left[T_{k-1}\left(1+2\,
 \frac{\la_i-\la_{p+1}}{\la_{p+1}-\la_n}\right)\right]^{-1}
 \tan\angle(\X,\Y),\quad i=1,\ldots,p,
\end{equation}
where the reciprocal Chebyshev term is the ratio
$|f(\la_{p+1})/f(\la_i)|$. Subsequently, since $f$ is a polynomial
of degree $k{-}1$ by \eqref{fsp}, the subspace $\Y'$ is a subset of
the block Krylov subspace $\K=\Y+A\Y+\cdots+A^{k-1}\Y$
(spanned by $k$ blocks) arising from the block Lanczos method.
Thus $\tan\angle(x_i,\K)\le\tan\angle(x_i,\Y')$ holds and
extends \eqref{blme} as a bound for $\tan\angle(x_i,\K)$.
The reciprocal Chebyshev term decreases rapidly with the degree $k$
of $\K$ provided that $\la_i$ is not close to $\la_{p+1}$.
This can be achieved by selecting an initial subspace of sufficiently high dimension. 

Bound \eqref{blme} demonstrates the quality of the shifted Chebyshev polynomial
\eqref{fsp} as a filter in \eqref{abi} amplifying the components of a signal $y\in\Y$
corresponding to the $p$ largest eigenvalues of $A$. For computing these eigenvalues,
one is also interested in bounds of the errors of the Ritz values of $A$ in $\Y'$,
arranged as $\eta'_1\ge\cdots\ge\eta'_p$. Some direct bounds
have been presented by Knyazev in \cite[Section 2]{ks}. In particular,
applying \cite[Theorem 2.5]{ks} to the specified function \eqref{fsp} yields the bound
\begin{equation}\label{blmervi}
 \frac{\la_i-\eta'_i}{\eta'_i-\la_n}\le\left[T_{k-1}\left(1+2\,
 \frac{\la_i-\la_{p+1}}{\la_{p+1}-\la_n}\right)\right]^{-2}
 \tan^2\angle(\X,\Y),\quad i=1,\ldots,p
\end{equation}
in the nontrivial case $\angle(\X,\Y)<\pi/2$. 
Although the convergence factor in \eqref{blmervi} is the squared value of that in \eqref{blme},
the derivation of \eqref{blmervi} is not simply based on \eqref{blme}. 
In addition, the relation $\Y'\subseteq\K$ leads to inequalities
$\eta'_i\le\psi_i$, $i=1,\ldots,p$ for the $p$ largest Ritz values
$\psi_1\ge\cdots\ge\psi_p$ of $A$ in $\K$ according to
the Courant-Fischer principles. This extends
\eqref{blmervi} to the block Lanczos method.
Bounds \eqref{blme} and \eqref{blmervi} are cluster robust, just as well as \eqref{bpme}.

In comparison to \cite[Theorems 5 and 6]{s1980} by Saad, bounds
\eqref{blme} and \eqref{blmervi} are much more accurate in the case where $\la_i$
belongs to an eigenvalue cluster; cf.~numerical examples
in \cite[Subsection 5.1]{z2018} concerning a slightly modified form of \eqref{blme}.
Indeed, the corresponding analysis in \cite{s1980}
can also be reformulated with respect to
the abstract block iteration \eqref{abi}. The specified function $f$ contains
linear factors of the form $(\alpha-\la_j)$ or $(\alpha-\psi_j)$,
$j=1,\ldots,i{\,-\,}1$ depending on eigenvalues or Ritz values
in order to construct auxiliary vectors which are orthogonal to the associated
eigenvectors or Ritz vectors. This leads to the terms
$\prod_{j=1}^{i-1}(\la_j-\la_n)/(\la_j-\la_i)$
and $\prod_{j=1}^{i-1}(\psi_j-\la_n)/(\psi_j-\la_i)$
called ``bulky'' factors in \cite{lz2015},
which could be very large for clustered eigenvalues
and cannot be influenced by the block size.

\subsection{Tuplewise analysis using majorization}

Majorization theory is a powerful tool that allows extend traditional bounds like
\eqref{bpmf}, \eqref{blme} and \eqref{blmervi} of individual convergence measures
to bounds of the corresponding tuples such as
\[[\tan\angle(x_c,\Y'),\ldots,\tan\angle(x_d,\Y')]
\quad\mbox{and}\quad
[\la_c-\eta'_c,\ldots,\la_d-\eta'_d]\]
for $1 \le c \le d \le p$.
Majorization theory has been first applied to analyzing
some subspace iterations and the block Lanczos method in \cite{ka10},
where convergence bounds involve tuples of convergence measures
and even the convergence factors can also be tuples.
Similar bounds of tuples with only scalar convergence factors
have been later independently derived in \cite{lz2015}.

Certain bounds from \cite{ka10}
are extendable to the convergence analysis of \eqref{abi}.
For instance, applying \cite[Theorem 2.3]{ka10}
yields a ``stationary'' bound in comparison to \eqref{blmervi}.
Let us consider the $p$ Ritz value errors $(\la_i-\eta'_i)/(\eta'_i-\la_n)$, $i=1,\ldots,p$,
and write them as a sorted $p$-tuple $\varepsilon$ whose components fulfill
$\varepsilon_1\ge\cdots\ge\varepsilon_p$. For bounding $\varepsilon$, we use
the arranged principal angles $\theta'_1\ge\cdots\ge\theta'_p$
between the subspaces $\X$ and $\Y'$.
Skipping the majorization notation, we get
\[\textstyle\sum_{j=1}^d \varepsilon_j\le\sum_{j=1}^d \tan^2\theta'_j,
 \quad d=1,\ldots,p.\]
An extension in terms of the arranged principal angles $\theta_1\ge\cdots\ge\theta_p$
between $\X$ and $\Y$ reads
\begin{equation}\label{blmervim}
 \textstyle\sum_{j=1}^d \varepsilon_j\le\sum_{j=1}^d \beta_j\tan^2\theta_j,
 \ \ d=1,\ldots,p \ \ \mbox{with} \ \ \displaystyle
 \beta_j=\left[T_{k-1}\left(1+2\,
 \frac{\la_{p+1-j}-\la_{p+1}}{\la_{p+1}-\la_n}\right)\right]^{-2}.
\end{equation}
Because $\theta_1=\angle(\X,\Y)$,
the special form \,$\varepsilon_1\le\beta_1\tan^2\theta_1$\,
of \eqref{blmervim} for $d=1$ can imply \eqref{blmervi} for $i=p$.
Moreover, summing up \eqref{blmervi} for all indices leads to the bound
\[\textstyle\sum_{j=1}^p \varepsilon_j
 \le\big(\sum_{j=1}^p \beta_j\big)\tan^2\theta_1.\]
The alternative bound
\[\textstyle\sum_{j=1}^p \varepsilon_j
 \le \beta_1\sum_{j=1}^p \tan^2\theta_j\]
is analogous to \cite[(8.11)]{lz2015}
where the convergence factor is also a scalar.
These two bounds can be improved by \eqref{blmervim} for $d{\,=\,}p$,
especially if the maximal component $\tan^2\theta_1$ or $\beta_1$ is dominant.

As a motivation of the present paper, it is remarkable that
the terms $\theta_j$ and $\beta_j$ in \eqref{blmervim}
are defined within $p$-tuples concerning
principal angles between two $p$-dimensional subspaces
and the first $p$ eigenvalues.
If \eqref{blmervim} is applied to some $d \ll p$,
a significant overestimation can occur since the utilized $\theta_j$ and $\beta_j$
are the $d$ largest components of the respective $p$-tuples.
Therefore, we expect to improve \eqref{blmervim} by using
certain terms which simply depend on $d$-dimensional subspaces
and the first $d$ eigenvalues.

\subsection{Aim and outline}

Deriving majorization bounds for the block Lanczos method
via the abstract block iteration \eqref{abi} is a project originated from \cite{ak08}.
Some previous results such as \eqref{blmervim} are included in Section 4
of the technical report \cite{ka08} but due to space limits removed
in the final version \cite{ka10}. The published results in \cite{ka10}
concern relations between Ritz value errors and principal angles,
but do not depend on \eqref{abi}. These results can also be extended to
Hermitian operators on infinite dimensional spaces; cf.~\cite[Subsection 2.8]{ka10}
and a recent application in a filtered subspace iteration \cite{ggo2020}. 
Moreover, some improvements or alternatives of the sine-based bound
\cite[(2.5)]{ka10} have been achieved in \cite{msz2020,msz2021}
and devoted to a partial confirmation of the conjecture \cite[(2.4)]{ka10}.
The remaining results in \cite{ka08}
actually enable a majorization-type generalization of bound \eqref{blmervi}.
Furthermore, the intermediate bound \cite[(4.1)]{ka08} is comparable with
\cite[Theorem 2.1]{dikm2018} concerning block Krylov subspace methods
for solving singular value problems. Another potential application
is the convergence analysis of a subspace iteration for polynomial
eigenvalue problems \cite{gpv2022}. The present paper extends \cite{ka08}
and makes further progress concerning certain biorthogonal vectors
which have been utilized implicitly for investigating the block power method
by Rutishauser; cf.~Lemma \ref{lm:auxvec}.

In the remaining part of this paper, we extend the existing convergence theory
of the abstract block iteration \eqref{abi} by majorization arguments.
In Section 2, we introduce basic settings and necessary definitions.
Section 3 provides main results after discussing some previous results
that provide some motivation for our investigation.
Applications to concrete methods
are formulated in Section 4. We mainly consider the block Lanczos method
and compare our new bounds with those from \cite{lz2015}.
Further applications are related to shift-and-invert eigensolvers
for certain discretized partial differential operators.
In Section 5, several numerical examples
illustrate the accuracy of the new bounds, followed by concluding remarks.
Section 6 (Appendix) presents some detailed proofs. 

\section{Preliminaries}

This section serves as preparation for our majorization-type convergence analysis of
block iterative eigensolvers. Subsection 2.1 contains
some basic settings which are consistently required in the analysis.
In Subsection 2.2, we introduce principal angles and their tangent definition
as a useful tool for deriving majorization-type bounds.
Subsection 2.3 provides tuple notations and majorization statements
for some singular value tuples. An overview of frequently used symbols is as follows.
\\
\mbox{\quad} Subsection 2.1: \ \
 matrices $A$, $I_d$, \ eigenpairs $(\la_i,x_i)$ of $A$,
 \ subspaces $\X$, $\Y$, $\Y'$.
\\
\mbox{\quad} Subsection 2.2: \ \
 principal angles, \ notation $\angle(\cdot,\cdot)$ for some special forms.
\\
\mbox{\quad} Subsection 2.3: \ \
 tuple $a^\downarrow$ with arranged components, \
 majorization relations $\prec_w$, $\prec$, \\
\mbox{\quad}\hspace{2.6cm}\ \
 singular value tuple $S(\cdot)$,
 eigenvalue tuple $\Lambda(\cdot)$,
 principal angle tuple $\Theta(\cdot,\cdot)$.

\subsection{Basic settings}

We consider a normal matrix $A\in\C^{n \times n}$ with the eigenvalues
$\la_1,\ldots,\la_n$ and the associated orthonormal eigenvectors $x_1,\ldots,x_n$.
Identity matrices are denoted by $I_d$
with the corresponding dimension $d\in\{1,\ldots,n\}$.
Some column vectors or ``tall'' block matrices are represented
with their components or submatrices separated by ``;''
within a row as in MATLAB/Octave notation.
The abstract block iteration $\Y'=f(A)\Y$ introduced in \eqref{abi}
serves to compute a moderate number of the first eigenvalues of $A$.

It is possible to frame our analysis with
a Hermitian operator $\mathcal{A}$ in a finite dimensional
inner product space $\mathcal{H}$ together with orthogonal projectors
of subspaces avoiding introducing any bases and thus shortening
the notation; cf.~Remark \ref{ov}. An extension to the infinite dimensional case
can be made as in \cite[Subsection 2.8]{ka10}.
We choose the matrix-based formulation since the new results in this paper
are devoted to the convergence theory of matrix eigensolvers
where the standard notation usually begins with matrices
(even if the implementation is actually matrix-free).
In addition, the resulting new bounds can easily be compared with
some existing ones for the block Lanczos method, e.g.,
those from \cite{s1980,lz2015}.

In practical applications, the dimension $p$ of the initial subspace $\Y\subset\C^n$
is larger than the number of the target eigenvalues, and much smaller than $n$. 
We define the invariant subspace $\X=\span\{x_1,\ldots,x_p\}$,
and assume that the eigenvalue sets $\{\la_1,\ldots,\la_p\}$ and
$\{\la_{p+1},\ldots,\la_n\}$ are disjoint so that $\X$ is unique.
Moreover, we assume $\angle(\X,\Y)<\pi/2$ in order to
exclude trivial terms such as $\tan\angle(\X,\Y)=\infty$.
We reuse the following assumption on the function $f$
made in \cite{k,ks},
\begin{equation}\label{ca}
 \max_{j\in\{p+1,\ldots,n\}}|f(\la_j)|<\min_{j\in\{1,\ldots,p\}}|f(\la_j)|,
\end{equation}
i.e., $f$ serves as a filter enlarging the components with $i\le p$
relative to those with $i>p$ and thus 
ensures that the corresponding convergence factors are smaller than $1$.

Our Ritz value analysis deals with a Hermitian $A$ and
the approximation of its largest eigenvalues. Therein we arrange the eigenvalues as
$\la_1\ge\cdots\ge\la_n$ and simply assume $\la_p>\la_{p+1}$
to ensure the uniqueness of $\X$. The analysis can easily be transformed
by the substitution $A \to -A$ to the consideration of the smallest eigenvalues.

\subsection{Principal angles}

We use principal angles to measure the distance from-to or between subspaces
of the invariant subspace $\X$ and the initial subspace $\Y$. In general, for two subspaces
$\U,\,\V\subset\C^n$ with $\dim\U\le\dim\V$ and their arbitrary
orthonormal basis matrices $U$ and $V$, the cosine values of the principal angles
from $\U$ to $\V$ are defined by the singular values of $V^HU$.
Therein the largest principal angle is the Euclidean angle $\angle(\U,\V)$
which can also be defined by
$\max_{u\in\U{\setminus}\{0\}}\min_{v\in\V{\setminus}\{0\}}\angle(u,v)$
with angles between nonzero vectors.

In the case $\dim\U=\dim\V$, one also says
the principal angles ``between $\U$ and $\V$''
and uses $\angle(\U,\V)=\angle(\V,\U)$
for denoting the largest principal angle. In the case $\dim\U=1\le\dim\V$,
there is only one principal angle. This coincides with the angle
$\angle(u,\V)=\min_{v\in\V{\setminus}\{0\}}\angle(u,v)$
for arbitrary $u\in\U{\setminus}\{0\}$. If also $\dim\V=1$,
the associated cosine and tangent values are actually
$|\cos\angle(u,v)|$ and $|\tan\angle(u,v)|$ for arbitrary $v\in\V{\setminus}\{0\}$.

For constructing majorization bounds in Section 3,
we utilize tangent values of principal angles
where an angle equal to $\pi/2$ would cause a trivial infinity bound.
Therefore we restrict the consideration to the case where $V^HU$ has full rank.
Then all singular values of $V^HU$ are nonzero so that
all principal angles from $\U$ to $\V$ are smaller than $\pi/2$.

\begin{remark}\label{ov}
A general consideration of tangent values of principal angles
within an operator-based formulation
has been presented in \cite[Theorem 3.1]{zk2013}. Alternatively, one can
apply the elegant form $P_{\VP}(P_{\V}P_{\U})^{\dag}$
with the corresponding orthogonal projectors; cf.~\cite[Theorem 4.1]{zk2013}.
Nevertheless, the following description
with a matrix product $\QVP^H\widetilde{U}(V^H\widetilde{U})^{\dag}$
is more appropriate for estimating error reductions with respect to tangent values
in a concise way without additional zero components.
\end{remark}

\begin{lemma}\label{lm:tangentform}
Consider the subspace $\V\subset\C^n$ and its orthogonal complement $\VP$
with their arbitrary orthonormal basis matrices
$V \in \C^{n \times t}$ and $\QVP \in \C^{n \times (n-t)}$.
Let $\widetilde{U} \in \C^{n \times s}$ with $s \le \min\{t,\,n{-}t\}$ be an arbitrary
(but not necessarily orthonormal) basis matrix of the subspace $\U\subset\C^n$
for which $V^H\widetilde{U}$ has full rank. Then the $s$ largest singular values of
the $(n{-}t){\times}t$ matrix $\QVP^H\widetilde{U}(V^H\widetilde{U})^{\dag}$
coincide with the tangent values of the principal angles from ${\U}$ to ${\V}$
where the symbol $\dag$ denotes the Moore-Penrose pseudoinverse. In particular,
$\|\QVP^H\widetilde{U}(V^H\widetilde{U})^{\dag}\|=\tan\angle(\U,\V)$
holds with the $2$-norm $\|\cdot\|$.
\end{lemma}
\noindent\textit{Proof.} $\to$ Subsection \ref{p:tangentform}.

The condition $s \le \min\{t,\,n{-}t\}$ in Lemma \ref{lm:tangentform}
is consistent with the analysis in further sections. Therein the subspace iterates
or their subsets are considered as $\U$, whereas an invariant subspace
of the same or higher dimension corresponds to $\V$, i.e., $s \le t$.
This leads to reasonable
characteristics of bounds such as the gap ratio $(\la_s-\la_{t+1})/(\la_{t+1}-\la_n)>0$
in a Chebyshev term. Moreover, $s \le n{-}t$ is naturally fulfilled in the context of
computing several eigenvalues of a large matrix or a matrix pair.

\subsection{Majorization}

We use majorization arguments for deriving bounds in terms of
certain tuples of real numbers concerning tangent values of principal angles
and Ritz value errors. For notational convenience we treat these tuples
as row vectors, e.g., $a=[a_1,\ldots,a_d]$
with the corresponding dimension $d$. In addition, we consider
a rearrangement of the components of \,$a$\, in (not strictly) decreasing order
and denote the resulting vector by $a^\downarrow$, i.e.,
$a_1^\downarrow \ge \cdots \ge a_d^\downarrow$.
Let $b$ be another tuple with dimension $d$. If the components
of $a^\downarrow$ and $b^\downarrow$ fulfill
\[\textstyle\sum_{i=1}^k a_i^\downarrow\le\sum_{i=1}^k b_i^\downarrow
 \quad\forall\ k\in\{1,\ldots,d\},\]
one says that $b$ weakly majorizes
(submajorizes) $a$. The weak (additive) majorization is denoted by $a \prec_w b$.
The corresponding strong majorization with the notation $a \prec b$
additionally means that the sum inequality for $k=d$ is actually an equality.
For nonnegative tuples $a$ and $b$ with dimension $d$,
the weak multiplicative majorization is defined by
\[\textstyle\prod_{i=1}^k a_i^\downarrow\le\prod_{i=1}^k b_i^\downarrow
 \quad\forall\ k\in\{1,\ldots,d\}\]
and denoted by $\log a \prec_w \log b$.
The strong version $\log a \prec \log b$ requires the additional condition
\,$\prod_{i=1}^d a_i^\downarrow=\prod_{i=1}^d b_i^\downarrow$.
Moreover, nonnegative tuples with different dimensions can be compared by
adding zeros to the shorter tuple.

Specific tuples in our analysis mostly consist of singular values
or eigenvalues of matrices. We generally consider a matrix
$B\in\C^{d \times s}$, and denote by $S(B)$
the tuple of arranged singular values of $B$ in decreasing order.
Multiple singular values are counted repeatedly so that $S(B)$
has dimension $\min\{d,s\}$. Similarly, we denote by $\Lambda(C)$ the tuple of
arranged eigenvalues of a Hermitian matrix
$C\in\C^{t \times t}$ in decreasing order. In particular, it holds that
$S(B)=\sqrt{\Lambda(B^HB)}$ for $d \ge s$ and
$S(B)=\sqrt{\Lambda(BB^H)}$ for $d \le s$
with componentwise square roots.

The arranged principal angles $\theta_1\ge\cdots\ge\theta_s$ from $\U$ to $\V$
give the angle tuple $\Theta(\U,\V)=[\theta_1,\ldots,\theta_s]$.
The corresponding tangent tuple
$\tan\Theta(\U,\V)=[\tan(\theta_1),\ldots,\tan(\theta_s)]$ is the leading $s$-subtuple
of $S\big(\QVP^H\widetilde{U}(V^H\widetilde{U})^{\dag}\big)$
according to Lemma \ref{lm:tangentform}.
In addition, by using an arbitrary orthonormal basis matrix $U$ of $\U$,
the cosine-type definition of principal angles shows
$\big(\cos\Theta(\U,\V)\big)^\downarrow=S(V^HU)$.

We complete this subsection by a basic argument for the majorization-type analysis.
\begin{lemma}\label{lm:major}
Consider the matrices $B_1\in\C^{d_1 \times d_2}$, $B_2\in\C^{d_2 \times d_3}$
and $B_3\in\C^{d_3 \times d_4}$. Let $S_t(B)$ be the leading $t$-subtuple
of the singular value tuple $S(B)$ for
$B\in\{B_1,\,B_2,\,B_3,\,B_1B_2B_3\}$ and $t\le\min\{d_1,d_2,d_3,d_4\}$.
Then it holds with componentwise multiplication, division
and power for $c\in\mathbb{N}$ of tuples, that
\,$S_t^c(B_1B_2B_3) \prec_w S_t^c(B_1)S_t^c(B_2)S_t^c(B_3)$\, and
\[S_t^c(B_1B_2B_3)/S_t^c(B_2) \prec_w S_t^c(B_1)S_t^c(B_3)
\quad\mbox{for}\quad S_t(B_2)>0.\]
\end{lemma}
\begin{proof}
Applying \cite[Theorem 4.4]{ka10} shows the expression
\[\log S(B_1B_2B_3)-\log S(B_2)\prec\log\big(S(B_1)S(B_3)\big)\]
which means that the singular value inequality
\[\textstyle\prod_{j=1}^k \sigma_{i_j}(B_1B_2B_3)
 \le \prod_{j=1}^k \big(\sigma_j(B_1)\sigma_j(B_3)\big)\sigma_{i_j}(B_2)\]
holds for each $k\in\{1,\ldots,d\}$ with $d=\max\{d_1,d_2,d_3,d_4\}$
and for each index set $\{i_1,\ldots,i_k\}\subseteq\{1,\ldots,d\}$
with $i_1<\cdots<i_k$. Therein zeros are occasionally added to
shorter singular value tuples to match the sizes,
and the equality is attained for $k=d$. 
By considering this singular value inequality to the $c$th power for
$k\in\{1,\ldots,t\}$ and $\{i_1,\ldots,i_k\}=\{1,\ldots,k\}$,
we get the weak multiplicative majorizations
\,$\log S_t^c(B_1B_2B_3)\prec_w\log \big(S_t^c(B_1)S_t^c(B_2)S_t^c(B_3)\big)$\, and
\[\log\big(S_t^c(B_1B_2B_3)/S_t^c(B_2)\big) \prec_w
 \log\big(S_t^c(B_1)S_t^c(B_3)\big)\quad\mbox{for}\quad S_t(B_2)>0.\]
These imply the weak (additive) majorizations in the claim of Lemma \ref{lm:major}
by applying the exponential function ``$\exp$'' which is convex and increasing;
see \cite[Proposition 4.B.2]{moa}.
\end{proof}
In the above proof, the composition ``$\exp\circ\log$'' is only applied to positive components
whereas zero components remain unchanged as they lead to a trivial case.
Alternatively, one can derive the underlying additive inequality directly by
the corresponding multiplicative inequality; cf.~\cite[Corollary 3.3.10]{hj2}
and \cite[Example II.3.5]{bhatia_book}.

\section{Main results}

We motivate our new analysis of the abstract block iteration \eqref{abi}
by introducing some auxiliary terms and arguments.
The basic ingredients are certain biorthogonal vectors
from $\X$ and $\Y$ inspired by the classical analysis
of the block power method by Rutishauser \cite[Theorem 2]{r1969}
(do not confuse these biorthogonal vectors with
those mentioned before \cite[Section 1]{r1969} which form two sets of iterates).
We note that the corresponding biorthogonality is represented
by a submatrix in \cite[(10)]{r1969} after some substitutions
and results in an asymptotic bound. Indeed, this biorthogonality
also enables a simpler proof of the nonasymptotic bound \cite[(14.11)]{p}
by Parlett. A similar argument in \cite[Lemma 4]{s1980} by Saad
is the starting point of an analysis of the block Lanczos method.

\begin{lemma}\label{lm:auxvec}
With the settings from Subsection 2.1, consider the orthonormal basis matrix
$X=[x_1,\ldots,x_p]$ of the invariant subspace $\X$ together with
an arbitrary orthonormal basis matrix $\widetilde{Y}\in\C^{n \times p}$ of $\Y$.
Then the $p{\times}p$ matrix $X^H\widetilde{Y}$ is invertible, and the vectors
\[y_i=\widetilde{Y}c_i, \quad i=1,\ldots,p
 \quad\mbox{with the columns $c_1,\ldots,c_p$ of $(X^H\widetilde{Y})^{-1}$}\]
are evidently linearly independent and form a basis matrix $Y=[y_1,\ldots,y_p]$.
Moreover, $x_1,\ldots,x_p$ and $y_1,\ldots,y_p$ are biorthogonal, i.e.,
\begin{equation}\label{biorth}
 x_i^Hy_j=\delta_{ij} \quad\mbox{for}\quad i,j\in\{1,\ldots,p\}.
\end{equation}
\end{lemma}
\begin{proof}
The smallest singular value of $X^H\widetilde{Y}$ coincides with $\cos\angle(\X,\Y)$
which is nonzero because $\angle(\X,\Y)<\pi/2$.
Thus $X^H\widetilde{Y}$ is invertible. The biorthogonality is verified by
\[x_i^Hy_j=(Xe_i)^H\widetilde{Y}c_j
 =e_i^HX^H\widetilde{Y}\big((X^H\widetilde{Y})^{-1}e_j\big)
 =e_i^He_j=\delta_{ij}\]
with the columns $e_1,\ldots,e_p$ of $I_p$.
\end{proof}

We introduce some single-angle bounds in Subsection 3.1 and present
the corresponding multi-angle majorization-type bounds
in Subsection 3.2. Some relatively long proofs are given in Section 6.

\subsection{Single-angle bounds}

In \cite{r1969,p}, the auxiliary vectors defined in Lemma \ref{lm:auxvec} are utilized
separately. The simple power method applied to them can be observed within
proper invariant subspaces according to the biorthogonality \eqref{biorth}. An essentially
analogous approach with the tangent description introduced in Lemma \ref{lm:tangentform}
leads to the following single-angle bound
which corresponds to an abstract form of \eqref{blme}.

\begin{lemma}\label{lm:blme}(cf.~\cite[Lemma 2.3.1]{k})
With the settings from Subsection 2.1 and Lemma \ref{lm:auxvec}, it holds that
\begin{equation}\label{blmea}
\begin{split}
 &\tan\angle(x_i,\Y') \le \sigma_i\tan\angle(x_i,y_i)
  \le \sigma_i\tan\angle(\X,\Y) \\[1ex]
 &\mbox{with}\quad
 \sigma_i=\frac{\max_{j\in\{p+1,\ldots,n\}}|f(\la_j)|}{|f(\la_i)|},
 \quad i=1,\ldots,p.
\end{split}
\end{equation}
\end{lemma}
\noindent\textit{Proof.} $\to$ Subsection \ref{p:blme}.

For the angle-dependent bound \eqref{blmervi} on 
the Ritz values $\eta'_1\ge\cdots\ge\eta'_p$ of $A$ in $\Y'$,
an abstract form reads
\begin{equation}\label{blmerviabi}
 \frac{\la_i-\eta'_i}{\eta'_i-\la_n}
 \le\frac{\max_{j\in\{p+1,\ldots,n\}}|f(\la_j)|^2}{\min_{j\in\{1,\ldots,i\}}|f(\la_j)|^2}
 \tan^2\angle(\X,\Y),\quad i=1,\ldots,p.
\end{equation}
Therein $\dim\Y'=p$ is ensured by the common assumption \eqref{ca};
see Lemma \ref{lm:blmervi} below.
For deriving \eqref{blmerviabi}, we can put several auxiliary vectors
from Lemma \ref{lm:auxvec} together
in order to construct proper auxiliary subspaces.
This improves the approach by Saad \cite[Theorem 6]{s1980}
which constructs auxiliary vectors with Ritz values and causes
less convenient terms in the bound. 

The first step of the derivation
of \eqref{blmerviabi} produces an angle bound for auxiliary subspaces.

\begin{lemma}\label{lm:blmervi}(cf.~\cite[Lemma 2.3.1]{k})
With the settings from Subsection 2.1 and Lemma \ref{lm:auxvec},
consider for an index $i\in\{1,\ldots,p\}$ the subspace
$\X_i$ with the orthonormal basis matrix $X_i=[x_1,\ldots,x_i]$ and
the subspace $\Y_i$ with the basis matrix $Y_i=[y_1,\ldots,y_i]$.
Then the matrix $Y'_i=f(A)Y_i$ has rank $i$ so that the subspace
$\Y'_i=\span\{Y'_i\}$ has dimension $i$. In particular,
the subspace iterate $\Y'$ of \eqref{abi} coincides with $\Y'_p$
and thus has dimension $p$. Moreover, it holds that
\begin{equation}\label{blmervia}
 \tan\angle(\X_i,\Y'_i)
 \le\frac{\max_{j\in\{p+1,\ldots,n\}}|f(\la_j)|}{\min_{j\in\{1,\ldots,i\}}|f(\la_j)|}
 \tan\angle(\X_i,\Y_i).
\end{equation}
\end{lemma}
\noindent\textit{Proof.} $\to$ Subsection \ref{p:blmervi}.

The next step of the derivation of \eqref{blmerviabi}
is to verify a direct relation between the convergence measures
$\tan\angle(\X_i,\Y'_i)$ and $(\la_i-\eta'_i)/(\eta'_i-\la_n)$.

\begin{lemma}\label{lm:blmervi1}(cf.~\cite[Lemmas 2.2.5 and 2.2.6]{k})
With the settings from Subsection 2.1, Lemma \ref{lm:auxvec}
and Lemma \ref{lm:blmervi}, consider the Ritz values
$\eta'_1\ge\cdots\ge\eta'_p$ of $A$ in $\Y'$. Then
\begin{equation}\label{blmervia1}
 \frac{\la_i-\eta'_i}{\eta'_i-\la_n}\le\tan^2\angle(\X_i,\Y'_i)
\end{equation}
holds for each $i\in\{1,\ldots,p\}$.
\end{lemma}
\noindent\textit{Proof.} $\to$ Subsection \ref{p:blmervi1}.

Based on \eqref{blmervia} and \eqref{blmervia1}, we get the bound
\begin{equation}\label{blmerviabi1}
 \frac{\la_i-\eta'_i}{\eta'_i-\la_n}
 \le\frac{\max_{j\in\{p+1,\ldots,n\}}|f(\la_j)|^2}{\min_{j\in\{1,\ldots,i\}}|f(\la_j)|^2}
 \tan^2\angle(\X_i,\Y_i),\quad i=1,\ldots,p.
\end{equation}
If the auxiliary subspaces $\X_i$ and $\Y_i$ need to be eliminated,
one can extend \eqref{blmerviabi1} as \eqref{blmerviabi}; see Subsection \ref{p:ext}.

\subsection{Multi-angle majorization-type bounds}

Our new bounds for the abstract block iteration \eqref{abi} are generalizations
of those introduced in Subsection 3.1.
We use partial sums of certain tuples
containing (squared) tangent values of principal angles instead of
the corresponding maxima.

We first generalize the single-angle bound \eqref{blmea}
to principal angles concerning a possible eigenvalue cluster
in the eigenvalue set $\{\la_1,\ldots,\la_p\}$.

\begin{theorem}\label{thm:blmemajor}
With the settings from Subsection 2.1 and Lemma \ref{lm:auxvec}, 
consider an index set $\tau=\{i_1,\ldots,i_t\}\subseteq\{1,\ldots,p\}$
with $i_1<\cdots<i_t$ as well as the subspaces
\[\X_{\tau}=\span\{x_{i_1},\ldots,x_{i_t}\}\quad\mbox{and}\quad
 \Y_{\tau}=\span\{y_{i_1},\ldots,y_{i_t}\}.\]
By using the notations from Subsection 2.3
for tuples and majorization, it holds that
\begin{equation}\label{blmeamajor}
 \tan\Theta(\X_{\tau},\Y') \prec_w
 \Phi_{\tau}\,\widehat{\Phi}_t\,
 \tan\Theta(\X_{\tau},\Y_{\tau})
\end{equation}
with $\Phi_{\tau}=\big[|f(\la_{i_1})|^{-1},\ldots,
 |f(\la_{i_t})|^{-1}\big]^{\downarrow},\ \
 \widehat{\Phi}=\big[|f(\la_{p+1})|,\ldots,|f(\la_n)|\big]^{\downarrow}
 \ \ \mbox{and} \ \ \widehat{\Phi}_t=\widehat{\Phi}(1{\,:\,}t)$.
\end{theorem}
\noindent\textit{Proof.} $\to$ Subsection \ref{p:blmemajor}.

Theorem \ref{thm:blmemajor} is especially suitable for the case where
the eigenvalues $\la_{i_1},\ldots,\la_{i_t}$ are
consecutive and clustered. The components of the
tuplewise convergence factor $\Phi_{\tau}\,\widehat{\Phi}_t$
can be bounded away from $1$ by selecting a proper $f$ such that
$\max_{j\in\{p+1,\ldots,n\}}|f(\la_j)|\ll\min_{j\in\{1,\ldots,p\}}|f(\la_j)|$.
An application to nonconsecutive eigenvalues is also feasible,
but not of practical interest.

\begin{remark}
A direct application of Lemmas \ref{lm:tangentform} and \ref{lm:major}
implies the more abstract bound
\begin{equation}\label{blmeamajorab}
 \tan\Theta(F\U,\V) \prec_w
 S_s\big((V^HFV)^{-1}\big)\,S_s\big(\QVP^HF\QVP\big)\,
 \tan\Theta(\U,\V)
\end{equation}
for a normal matrix $F\in\C^{n \times n}$
provided that $\V$ and $\VP$ are invariant subspaces of $F^H$
and that $V^HFV$ is invertible;
see Subsection \ref{p:blmeamajorab} for the derivation.

We note that bound \eqref{blmeamajorab} can directly lead to
bound \eqref{blmeamajor} only in the special case
$\tau=\{1,\ldots,p\}$ because the biorthogonality \eqref{biorth} concerning
the index set $\tau$ cannot be added to \eqref{blmeamajorab} afterwards.
Alternatively, we can formulate a restricted form of \eqref{blmeamajorab}
with respect to an invariant subspace of $A$ associated with the eigenvalues
$\la_{i_1},\ldots,\la_{i_t},\la_{p+1},\ldots,\la_n$.
Nevertheless, we prefer the easily understandable formulation
in Theorem \ref{thm:blmemajor}.
\end{remark}

In order to eliminate $\Y_{\tau}$,
bound \eqref{blmeamajor} can be modified as
\begin{equation}\label{blmeamajor1a}
 \tan\Theta(\X_{\tau},\Y') \prec_w
 \Phi_{\tau}\,\widehat{\Phi}_t\,
 \tan\Theta_t(\X,\Y)
\end{equation}
with the leading $t$-subtuple $\Theta_t(\X,\Y)$ of $\Theta(\X,\Y)$.
Therein Lemma \ref{lm:tangentform} shows
\[\tan\Theta(\X_{\tau},\Y_{\tau})=S(\widetilde{G})\quad\mbox{for}\quad
 \widetilde{G}=[x_{p+1}^HY_{\tau};\,\ldots;\,x_n^HY_{\tau}]\in\C^{(n-p) \times t}\]
\[\mbox{and}\quad\tan\Theta(\X,\Y)=S(\widehat{G})\quad\mbox{for}\quad
 \widehat{G}=[x_{p+1}^HY;\,\ldots;\,x_n^HY]\in\C^{(n-p) \times p}\]
as in Subsection \ref{p:ext}.
Moreover, the Courant-Fischer principles ensure the tuple inequality
\[S(\widetilde{G})=\sqrt{\Lambda\big(\widetilde{G}^H\widetilde{G}\big)}
 =\sqrt{\Lambda\big(E_{\tau}^H\widehat{G}^H\widehat{G}E_{\tau}\big)}
 \le\sqrt{\Lambda_t\big(\widehat{G}^H\widehat{G}\big)}=S_t(\widehat{G}),\]
where $E_{\tau}$ is a $p{\times}t$ orthonormal matrix
consisting of the columns of $I_p$
with indices in $\tau$, and $\Lambda_t$ or $S_t$ denotes
the leading $t$-subtuple of the corresponding eigenvalue tuple
or singular value tuple. This leads to
\,$\Phi_{\tau}\widehat{\Phi}_t\tan\Theta(\X_{\tau},\Y_{\tau})$
$\le\Phi_{\tau}\widehat{\Phi}_t\tan\Theta_t(\X,\Y)$\,
and verifies \eqref{blmeamajor1a}.

For deriving a majorization version of the angle-dependent bound
\eqref{blmerviabi}, we proceed in two steps as in Subsection 3.1.
The first step gives an intermediate multi-angle bound
for auxiliary subspaces, and the second step enables a combination
between two convergence measures.

We formulate the first step as a generalization of Lemma \ref{lm:blmervi}
in square form.

\begin{lemma}\label{lm:blmervimajor}
(majorization update of \cite[Lemma 2.3.1]{k})
With the settings from Subsection 2.1, Lemma \ref{lm:auxvec}
and Lemma \ref{lm:blmervi}, it holds that
\begin{equation}\label{blmerviamajor}
 \tan^2\Theta(\X_i,\Y'_i) \prec_w
 \Phi_i^2\,\widehat{\Phi}_i^2\,
 \tan^2\Theta(\X_i,\Y_i)
\end{equation}
with $\Phi_i=\big[|f(\la_1)|^{-1},\ldots,
 |f(\la_i)|^{-1}\big]^{\downarrow},\ \
 \widehat{\Phi}=\big[|f(\la_{p+1})|,\ldots,|f(\la_n)|\big]^{\downarrow}
 \ \ \mbox{and} \ \ \widehat{\Phi}_i=\widehat{\Phi}(1{\,:\,}i)$.
\end{lemma}
\begin{proof}
Bound \eqref{blmerviamajor} is derived by adapting
the proof of Theorem \ref{thm:blmemajor}
(see Subsection 6.5) to $\tau=\{1,\ldots,i\}$ and $t=i$.
Most arguments before the intermediate bound \eqref{blmeamajor1}
are used in the same way. Only the application of Lemma \ref{lm:major}
is slightly different, namely, by using $c=2$ instead of $c=1$.
\end{proof}

The second step corresponds to a generalization of Lemma \ref{lm:blmervi1}.

\begin{lemma}\label{lm:blmervimajor1}
(majorization update of \cite[Lemmas 2.2.5 and 2.2.6]{k})
With the settings from Subsection 2.1, Lemma \ref{lm:auxvec},
Lemma \ref{lm:blmervi} and Lemma \ref{lm:blmervi1},
\begin{equation}\label{blmerviamajor1}
 \left[\frac{\la_1-\eta'_1}{\eta'_1-\la_n},\,\ldots,\,
 \frac{\la_i-\eta'_i}{\eta'_i-\la_n}\right]
 \prec_w \tan^2\Theta(\X_i,\Y'_i)
\end{equation}
holds for each $i\in\{1,\ldots,p\}$.
\end{lemma}
\noindent\textit{Proof.} $\to$ Subsection \ref{p:blmervimajor1}.

Now the majorization version of \eqref{blmerviabi} can be shown.
\begin{theorem}\label{thm:blmervimajor}
With the settings from Subsection 2.1, Lemma \ref{lm:auxvec},
Lemma \ref{lm:blmervi} and Lemma \ref{lm:blmervi1}, it holds that
\[\left[\frac{\la_1-\eta'_1}{\eta'_1-\la_n},\,\ldots,\,
 \frac{\la_i-\eta'_i}{\eta'_i-\la_n}\right]
 \prec_w\Phi_i^2\,\widehat{\Phi}_i^2\,\tan^2\Theta(\X_i,\Y_i)
 \le\Phi_i^2\,\widehat{\Phi}_i^2\,\tan^2\Theta_i(\X,\Y)\]
with $\Phi_i=\big[|f(\la_1)|^{-1},\ldots,
 |f(\la_i)|^{-1}\big]^{\downarrow},\ \
 \widehat{\Phi}=\big[|f(\la_{p+1})|,\ldots,|f(\la_n)|\big]^{\downarrow},
 \ \ \widehat{\Phi}_i=\widehat{\Phi}(1{\,:\,}i)$
and the leading $i$-subtuple $\Theta_i(\X,\Y)$ of $\Theta(\X,\Y)$.
\end{theorem}
\begin{proof}
The weak majorization is verified by combining \eqref{blmerviamajor}
with \eqref{blmerviamajor1}. The tuple inequality is shown analogously
to the derivation of \eqref{blmeamajor1a}.
\end{proof}

Theorem \ref{thm:blmervimajor} additionally provides
a majorization bound with $\Phi_i^2\,\widehat{\Phi}_i^2\,\tan^2\Theta_i(\X,\Y)$
which does not depend on auxiliary subspaces.

\begin{remark}
A more abstract bound concerning Lemma \ref{lm:tangentform} 
in the case $\dim\U=s=t=\dim\V$ can be derived by combining
\eqref{blmeamajorab} with an analogue of \eqref{blmerviamajor1}.
Therein $\V$ is assumed to be an invariant subspace of $A$
associated with the $t$ largest eigenvalues $\la_1\ge\cdots\ge\la_t$.
Then it holds that
\begin{equation}\label{blmerviamajorab}
 \left[\frac{\la_1-\psi_1}{\psi_1-\psi},\,\ldots,\,
 \frac{\la_t-\psi_t}{\psi_t-\psi}\right]
 \prec_w S^2\big((V^HFV)^{-1}\big)\,S_t^2\big(\QVP^HF\QVP\big)\,
 \tan^2\Theta(\U,\V)
\end{equation}
for the Ritz values $\psi_1\ge\cdots\ge\psi_t$ of $A$ in $F\U$
and the smallest Ritz value $\psi$ of $A$ in $F\U{\,+\,}\V$.
Proving Theorem \ref{thm:blmervimajor} based on \eqref{blmerviamajorab}
is only feasible for $i=p$ since the terms $|f(\la_{i+1})|,\ldots,|f(\la_p)|$
cannot easily be dropped afterwards.

Furthermore, we can use
\cite[Theorem 2.1]{ka10} for formulating an alternative bound
under the weaker assumption that
$\V$ is an invariant subspace of $A$ but not necessarily
associated with the $t$ largest eigenvalues. More precisely,
we reformulate \cite[(2.2)]{ka10} as
\[\zeta^{-1}\left| \Lambda\left(V^HAV\right) - \Lambda\left(W^HAW\right) \right|
 \prec_w \sin^2\Theta(F\U,\V).\]
with an orthonormal basis matrix $W$ of $F\U$ and the spread $\zeta$
of the Ritz value set of $A$ in $F\U{\,+\,}\V$.
Subsequently, we apply the convex and increasing function
$\alpha/(1-\alpha)$ defined for $\alpha\in[0,1)$ to the tuples
in the above majorization-type bound. Then
\[\frac{\left|\Lambda\left(V^HAV\right)-\Lambda\left(W^HAW\right)\right|}
 {\,[\zeta,\ldots,\zeta]-
 \left|\Lambda\left(V^HAV\right)-\Lambda\left(W^HAW\right)\right|\,}
 \prec_w \tan^2\Theta(F\U,\V).\]
Combining this with \eqref{blmeamajorab} results in
\begin{equation}\label{blmerviamajorab1}
 \frac{\left|\Lambda\left(V^HAV\right)-\Lambda\left(W^HAW\right)\right|}
 {\,[\zeta,\ldots,\zeta]-
 \left|\Lambda\left(V^HAV\right)-\Lambda\left(W^HAW\right)\right|\,}
 \prec_w S^2\big((V^HFV)^{-1}\big)S_t^2\big(\QVP^HF\QVP\big)
 \tan^2\Theta(\U,\V).
\end{equation}
The weaker assumption on $\V$ enables the application of \eqref{blmerviamajorab1}
to interior eigenvalues. However, it is somewhat challenging to select
an optimal filter $F$ without knowing good approximations
to undesired eigenvalues for polynomial filters
and/or to desired eigenvalues for rational filters.
The well-known Chebyshev terms for this type of bounds actually require that
$\V$ corresponds to the contiguous set of the $t$ largest eigenvalues.
\end{remark}

\section{Applications to specific filters}

We apply our new bounds presented in Subsection 3.2 to the convergence analysis
of the block Lanczos method and its inverted version.
The construction of the specified bounds essentially consists of
specifying the filter $f$ in the abstract block iteration \eqref{abi}
and adapting intermediate bounds to explicit convergence measures
based on subspace inclusions.
The practical choices of $f$ are normally limited to polynomials and rational functions.
The polynomial version of \eqref{abi} is easier to implement, since
its essential part is the multiplication of vectors or basis matrices by $A$.
For implementing the rational version of \eqref{abi},
one requires linear system solvers which make a single step
more expensive with respect to computational time and storage requirements.
Nevertheless, the total computational time can significantly be reduced.

In Subsection 4.1, we specify the filter $f$ as various shifted Chebyshev polynomials
so that \eqref{abi} simulates certain underlying iterations with almost optimal
convergence rates within the block Lanczos method.
This is comparable with the standard approach
from \cite{s1980,p,k} and the recent approach from \cite{lz2015}.
A remarkable feature of our majorization-type bounds
is that the convergence factors are tuples instead of scalars
as in the existing bounds.
In Subsection 4.2, we introduce an adaption of the main results
to some shift-and-invert eigensolvers.

\subsection{Block Lanczos method}

With the settings from Subsection 2.1,
the block Lanczos method constructs
block Krylov subspaces of the form $\K=\Y+A\Y+\cdots+A^{k-1}\Y$.
Thus the abstract block iteration \eqref{abi}, with an arbitrary real polynomial
of degree $k{-}1$ as $f$, produces a subset $\Y'$ of $\K$.
The inclusion $\Y'\subseteq\K$ leads to simple
inequalities for principal angles and Ritz values.

In order to specify $f$ as a reasonable polynomial, we consider
the multi-angle bound \eqref{blmeamajor} as an example.
Ideally, we want to construct $f$ in an optimal way, namely, minimizing the
tuplewise convergence factor $\Phi_{\tau}\,\widehat{\Phi}_t$ defined by
\[\Phi_{\tau}=\big[|f(\la_{i_1})|^{-1},\ldots,
 |f(\la_{i_t})|^{-1}\big]^{\downarrow},\ \
 \widehat{\Phi}=\big[|f(\la_{p+1})|,\ldots,|f(\la_n)|\big]^{\downarrow}
 \ \ \mbox{and} \ \ \widehat{\Phi}_t=\widehat{\Phi}(1{\,:\,}t)\]
with respect to the included eigenvalues. Following the standard approach
from \cite{s1980,p,k}, we simplify this minimization problem
with respect to the eigenvalue interval $[\la_n,\la_{p+1}]$.
Therein we begin with the tuple inequality
\begin{equation}\label{minbound}
 \Phi_{\tau}\,\widehat{\Phi}_t\le\Phi_{\tau}\,\varphi
 =\big[|f(\la_{i_1})|^{-1}\varphi,\ldots,|f(\la_{i_t})|^{-1}\varphi\big]^{\downarrow}
 \quad\mbox{with}\quad
 \varphi=\max_{\la\in[\la_n,\la_{p+1}]}|f(\la)|
\end{equation}
and then determine an $f$ minimizing the upper bound $\Phi_{\tau}\,\varphi$.
Fortunately, a shifted Chebyshev polynomial simultaneously minimizes 
all the components of $\Phi_{\tau}\,\varphi$
as shown in the following lemma inspired by \cite[Lemma 2.4.1]{k}.

\begin{lemma}\label{lm:scp}
Every component of the tuple $\Phi_{\tau}\,\varphi$ in \eqref{minbound}
is minimized in the class of real polynomials of degree $k{-}1$ at
\begin{equation}\label{scp}
 f(\alpha)=T_{k-1}\left(\frac{2\alpha-\la_{p+1}-\la_n}{\la_{p+1}-\la_n}\right)
 =T_{k-1}\left(1+2\,\frac{\alpha-\la_{p+1}}{\la_{p+1}-\la_n}\right)
\end{equation}
where $T_l$ with $l\in\mathbb{N}$ denotes the Chebyshev polynomials
(of the first kind). With this choice of $f$, we have
\,$\Phi_{\tau}\,\varphi=[\sigma_{i_t},\ldots,\sigma_{i_1}]$\, where
\begin{equation}\label{scpbound}
\begin{split}
 &\quad\sigma_j=\left[T_{k-1}\left(\frac{1+\xi_j}{1-\xi_j}\right)\right]^{-1}
 =\left[T_{k-1}(1+2\gamma_j)\right]^{-1}\\[1ex]
 &\mbox{with}\quad
 \xi_j=\frac{\la_j-\lambda_{p+1}}{\la_j-\la_n}\quad\mbox{and}\quad
 \gamma_j=\frac{\la_j-\la_{p+1}}{\la_{p+1}-\la_n}. 
\end{split}
\end{equation}
\end{lemma}
\begin{proof}
The absolute value of the Chebyshev polynomial $T_l$
is bounded above by $1$ on the interval $[-1,1]$
and exceeds $1$ outside of this interval
where the growth is faster in comparison to any other polynomial
of degree $l$ whose absolute value is also bounded above by $1$ on $[-1,1]$.
In our context, after mapping the interval $[-1,1]$ linearly to
the interval $[\la_n,\la_{p+1}]$, the above well-known property indicates that
function \eqref{scp} solves the minimization problem
described in the lemma. The verification of
$\Phi_{\tau}\,\varphi=[\sigma_{i_t},\ldots,\sigma_{i_1}]$
is straightforward where the arrangement of the components
is based on the monotonicity of $T_l$ on $(1,\infty)$.
\end{proof}

Several implementations of subspace iterations directly using
Chebyshev polynomials are available; see, e.g., \cite{p}
for the three-term recurrence and \cite{ks}
for the two-term recurrence iterative formulas.
Therein reasonable quality bounds of $\la_{p+1}$ and $\la_n$
are normally required.

The convergence theory
of subspace iterations using Chebyshev polynomials
is a standard way to derive the convergence rate bounds
for the block Lanczos method. We note that
the $\xi_j$-notation and the $\gamma_j$-notation
of the Chebyshev term \eqref{scpbound}
are suggested in \cite{ks} and \cite{s1980}, respectively.
In the $\gamma_j$-notation, one can easily see that
the bound decreases if the so-called gap ratio $\gamma_j$ increases.
The sharpness of scalar Chebyshev bounds for the single-vector version
of the Lanczos method has been discussed in \cite{l2010,zn2017}.
Some more accurate bounds are constructed by interpolating polynomials
in \cite{zn2017} and adapted to the block Lanczos method in \cite{z2018}.

We continue with the reformulation of majorization bounds.
Based on Lemma \ref{lm:scp}, we specify the bounds from Subsection 3.2
for the block Lanczos method as follows.

\begin{theorem}\label{thm:bl}
With the settings from Subsection 2.1 and Lemma \ref{lm:auxvec}, 
consider an index set $\tau=\{i_1,\ldots,i_t\}\subseteq\{1,\ldots,p\}$
with $i_1<\cdots<i_t$ as well as the subspaces
\[\X_{\tau}=\span\{x_{i_1},\ldots,x_{i_t}\},\quad
 \Y_{\tau}=\span\{y_{i_1},\ldots,y_{i_t}\}\quad\mbox{and}\quad
 \K=\Y+A\Y+\cdots+A^{k-1}\Y.\]
By using the notations from Subsection 2.3
for tuples and majorization, and the parameter
definition \eqref{scpbound}, it holds that
\begin{equation}\label{blmeamajorbl}
 \tan\Theta(\X_{\tau},\K)
 \prec_w[\sigma_{i_t},\ldots,\sigma_{i_1}]\,\tan\Theta(\X_{\tau},\Y_{\tau})
 \le[\sigma_{i_t},\ldots,\sigma_{i_1}]\,\tan\Theta_t(\X,\Y)
\end{equation}
with the leading $t$-subtuple $\Theta_t(\X,\Y)$ of $\Theta(\X,\Y)$.

Next, consider the $p$ largest Ritz values $\psi_1\ge\cdots\ge\psi_p$
of $A$ in $\K$ together with the subspaces
$\X_i=\span\{x_1,\ldots,x_i\}$ and $\Y_i=\span\{y_1,\ldots,y_i\}$
for an index $i\in\{1,\ldots,p\}$. Then
\begin{equation}\label{blmervimajorbl}
 \left[\frac{\la_1-\psi_1}{\psi_1-\la_n},\,\ldots,\,
 \frac{\la_i-\psi_i}{\psi_i-\la_n}\right]
 \prec_w [\sigma_i^2,\ldots,\sigma_1^2]\,\tan^2\Theta(\X_i,\Y_i)
 \le[\sigma_i^2,\ldots,\sigma_1^2]\,\tan^2\Theta_i(\X,\Y)
\end{equation}
holds with the leading $i$-subtuple $\Theta_i(\X,\Y)$ of $\Theta(\X,\Y)$.
\end{theorem}
\begin{proof}
For any real polynomial $f$ of degree $k{-}1$
we have $\Y'=f(A)\Y\subseteq\K$ and thus
\[\tan\Theta(\X_{\tau},\K)\le\tan\Theta(\X_{\tau},\Y'),\quad
 \la_j\ge\psi_j\ge\eta'_j \ \ \forall\ j\in\{1,\ldots,p\},\]
where the tangent tuple inequality can be proved analogously
to the end of the proof of Theorem \ref{thm:blmemajor},
and the Ritz value inequalities are ensured by the Courant-Fischer principles.
Then the left-hand sides of the bounds from
Theorem \ref{thm:blmemajor}, its supplement \eqref{blmeamajor1a}
and Theorem \ref{thm:blmervimajor}
can be extended as the left-hand sides of the corresponding specified bounds.
The extension of the right-hand sides is justified by Lemma \ref{lm:scp}.
\end{proof}

\begin{remark}
The multi-angle bound \eqref{blmeamajorbl} is a majorization-type
generalization of the single-angle bound \cite[(15)]{z2018} which
corresponds to \eqref{blme} and improves \cite[(3.4)]{s1980}. 
A similar multi-angle bound reads
\begin{equation}\label{blmeamajorblui}
 \textstyle\sum_{j=1}^l\tan\theta_j(\X_{\tau},\K)
 \le\sum_{j=1}^l\sigma_{i_t}\tan\theta_j(\X_{\tau},\Y_{\tau}),\quad
 l=1,\ldots,t
\end{equation}
by applying \cite[Theorem 8.1]{lz2015} where the underlying
unitarily invariant norm is specified as the Ky~Fan~$l$-norm
for $l=1,\ldots,t$, i.e., the sum of the $l$ largest singular values;
cf.~\cite[(IV.33)]{bhatia_book}.
In comparison to \eqref{blmeamajorblui},
the new bound \eqref{blmeamajorbl} is more accurate
by considering that $\sigma_{i_t}$ is the maximal component
of the tuplewise convergence factor $[\sigma_{i_t},\ldots,\sigma_{i_1}]$.
\end{remark}

\begin{remark}\label{bks}
We can modify \eqref{blmeamajorbl} by restricting its derivation
to an invariant subspace which is orthogonal to the eigenvectors
associated with a number of the largest eigenvalues
and occasionally also the smallest eigenvalues;
cf.~\cite{s1980} and \cite[Section 12.5]{p}. Then additional terms like
$\prod_{j=1}^{i-1}(\la_j-\la_n)/(\la_j-\la_i)$ occur in the bound
so that this modification is only meaningful for well-separated eigenvalues.
If the target eigenvalues are clustered and the initial subspace in
the block Lanczos method is relatively small, we can interpret
\eqref{blmeamajorbl} in another way, namely, let $\Y$ be
a block Krylov subspace
$\widetilde{\K}=\widetilde{\Y}+A\widetilde{\Y}+\cdots+A^{l-1}\widetilde{\Y}$
(therein $\dim\widetilde{\K} \le l\,\dim\widetilde{\Y}$ holds,
not necessarily with equality),
then $\K=\Y+A\Y+\cdots+A^{k-1}\Y$ is related to
a block Krylov subspace with the initial subspace $\widetilde{\Y}$
and degree $k{+}l{-}1$. This improves the applicability of
\eqref{blmeamajorbl} concerning small initial subspaces;
cf.~\cite[Section 8]{lz2015}.
\end{remark}

In addition, the Ritz vector bounds \cite[(4.10) and (8.8)]{lz2015}
are generalizations of \cite[(2.15)]{s1980} related to subspaces spanned
by Ritz vectors with a unitarily invariant norm.
We note that a generalization with a standard operator norm
is presented in \cite[Theorem 4.3]{k97}.
This bound type requires Ritz values.
A drawback is that the bound can decrease if some estimated Ritz values are
utilized instead of the corresponding exact Ritz values so that a combination with
Ritz value bounds is not meaningful. For this reason, we recommend
bound \cite[(2.7)]{ks}. Moreover, direct and concise Ritz vector bounds
are perhaps only known for the single-vector version of the Lanczos method;
see \cite[(3.3)]{nz2016}. Improving Ritz vector bounds
is a potential topic in our future research.

We now turn our attention to the angle-dependent
Ritz value bound \eqref{blmervimajorbl}.
Similarly to Remark \ref{mod}, a slightly sharper
form of \eqref{blmervimajorbl} with the smallest Ritz value
of $A$ in $\X{\,+\,}\K$ instead of $\la_n$ can be constructed
by a restricted analysis in $\X{\,+\,}\K$.
Indeed, also the parameter definition \eqref{scpbound}
can be modified for \eqref{blmervimajorbl} 
by using certain Ritz values of $A$ in $\X{\,+\,}\K$
instead of $\la_{p+1}$ and $\la_n$. This is enabled by the following fact:
With an arbitrary orthonormal basis matrix $V$ of $\X{\,+\,}\K$, the relation
\[\begin{split}
 V^HA^j\Y&=V^HA(A^{j-1}\Y)=V^HA(VV^H)(A^{j-1}\Y)\\
 &=(V^HAV)V^HA^{j-1}\Y=\cdots=(V^HAV)^j(V^H\Y)
\end{split}\]
holds for each $j\in\{1,\ldots,k{-}1\}$
by using the orthogonal projector $VV^H$. Thus 
\[\begin{split}
 V^H\K&=V^H\Y+V^HA\Y+\cdots+V^HA^{k-1}\Y\\
 &=V^H\Y+(V^HAV)(V^H\Y)+\cdots+(V^HAV)^{k-1}(V^H\Y)
\end{split}\]
is a block Krylov subspace with respect to $V^HAV$,
and analogous bounds can be achieved with eigenvalues
of $V^HAV$, i.e., Ritz values of $A$ in $\X{\,+\,}\K$.
This fact is a direct generalization of that
for Krylov subspaces; cf.~\cite[pp.~36]{k} and \cite[Lemma 3.3]{nz2016}.

\begin{remark}
Bound \eqref{blmervimajorbl} generalizes \cite[(2.20)]{ks}
which can be reformulated as \eqref{blmervi} or \cite[(19)]{z2018}
and improves \cite[(3.10)]{s1980}.
A similar bound based on \cite[Theorem 8.2]{lz2015} reads
\begin{equation}\label{blmervimajorblui}
 \textstyle\sum_{j=1}^l\dfrac{\la_j-\psi_j}{\la_1-\la_n}
 \le\sum_{j=1}^l\sigma_i^2\tan^2\theta_j(\X_i,\Y_i),\quad
  l=1,\ldots,i.
\end{equation}
Comparing the new bound \eqref{blmervimajorbl} with \eqref{blmervimajorblui}
indicates a twofold improvement
because $(\la_1-\la_n)^{-1}\le(\psi_j-\la_n)^{-1}$ and
$\sigma_i^2\ge\sigma_j^2\ \ \forall\ j\in\{1,\ldots,i\}$.

In addition, setting auxiliary vectors orthogonal to the Ritz vectors
associated with a number of the largest Ritz values leads to
an alternative bound with terms like
$\prod_{j=1}^{i-1}(\psi_j-\la_n)/(\psi_j-\la_i)$
which are suboptimal for clustered eigenvalues.
Furthermore, the interpretation of \eqref{blmeamajorbl} in Remark \ref{bks}
concerning small initial subspaces also fits \eqref{blmervimajorbl}.
\end{remark}

\subsection{Shift-and-invert eigensolvers}

We consider a generalized eigenvalue problem $Lv = \alpha Sv$ arising from
the finite element discretization of an operator eigenvalue problem.
Therein $L$ and $S$ are $n{\times}n$ Hermitian matrices,
and $S$ is positive definite. This formally includes the case of
the finite difference discretization with $S=I_n$.
Usually one only needs to compute a moderate number of eigenvalues.

By using a proper shift $\beta$, the shifted matrix
$L_{\beta}=L-\beta S$ is invertible, and computing
eigenvalues of $(L,S)$ close to $\beta$ corresponds to
computing extremal eigenvalues of $(S,L_{\beta})$.
The latter problem can be reformulated (implcitly) as computing
the largest eigenvalues of the Hermitian matrix pair 
\[(\widetilde{L},M)
 \quad\mbox{with}\quad
 \widetilde{L}=\pm L_{\beta}
 \quad\mbox{and}\quad
 M=L_{\beta}S^{-1}L_{\beta}\]
where $M$ is positive definite. This is equivalent to 
computing the largest eigenvalues of
\[A=M^{-1/2}\widetilde{L}M^{-1/2}.\]
A reverse transformation toward the original problem
can be used to construct shift-and-invert versions
of various iterative eigensolvers.

For instance, a block Krylov subspace $\K=\Y+A\Y+\cdots+A^{k-1}\Y$
can be transformed as
\[\widehat{\K}=\Z+M^{-1}\widetilde{L}\Z+\cdots+(M^{-1}\widetilde{L})^{k-1}\Z
 \quad\mbox{with}\quad\widehat{\K}=M^{-1/2}\K
 \quad\mbox{and}\quad\Z=M^{-1/2}\Y.\]
Such block Krylov subspaces with respect to $M^{-1}\widetilde{L}$
correspond to a shift-and-invert version of the block Lanczos method.
For an implementation with practical construction of $\widehat{\K}$,
we can solve linear systems of the form $Mw=r$
for certain Ritz vector residuals $r$ similarly to the generalized Davidson method.
The linear system $Mw=r$ is actually $L_{\beta}(S^{-1}L_{\beta}w)=r$
and can thus be solved as two successive linear systems for
the shifted matrix $L_{\beta}$.

For the convergence analysis, we can reformulate Lemma \ref{lm:scp}
and Theorem \ref{thm:bl} based on the above substitutions.
Therein the parameter definition \eqref{scpbound}
and the Ritz value measures are reformulated
by $\la=\pm(\alpha-\beta)^{-1}$ for eigenvalues.
The transformation of angle terms 
is related to the inner product induced by $M$.

\section{Numerical examples}

We discuss the accuracy of our new results with several numerical examples.
Examples 1 and 2 deal with the block Lanczos method applied to real diagonal matrices
following a classical example from \cite{s1980}. The bounds from Theorem \ref{thm:bl}
and their counterparts from \cite{lz2015} can directly be applied and demonstrated
with a comparative illustration.
The associated MATLAB/Octave codes (including validity checks
for Theorems \ref{thm:blmemajor} and \ref{thm:blmervimajor}) are available on
\begin{center}
\verb|https://github.com/lobpcg/MAJORIZATION_TYPE_CLUSTER_ROBUST_BOUNDS|\\
\verb|_FOR_BLOCK_FILTERS_AND_EIGENSOLVERS|
\end{center}

\subsection{Example 1}
\begin{figure}[htbp]
\begin{center}
\includegraphics[width=0.9\textwidth]{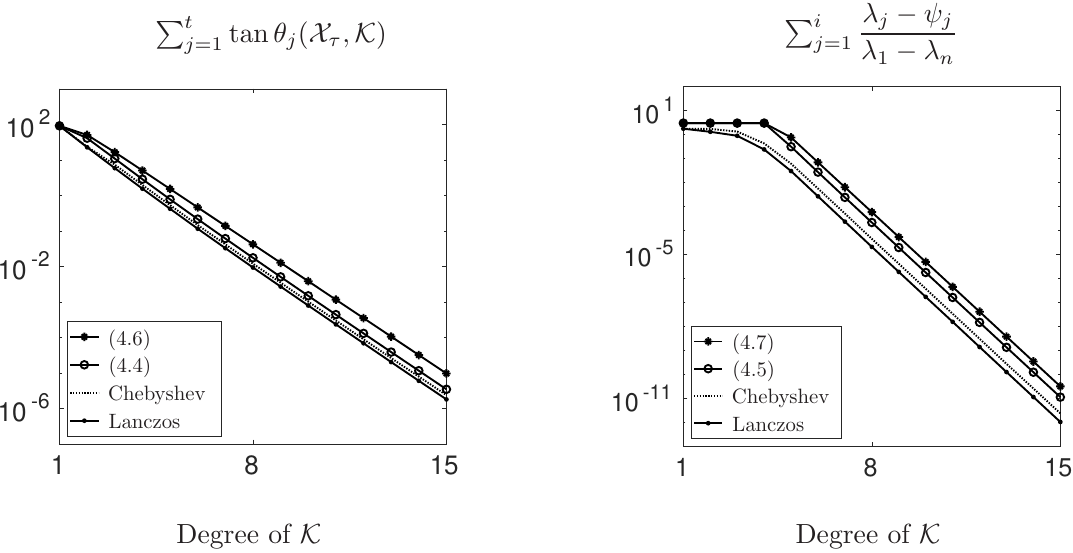}\quad
\end{center}
\par\vskip -2ex
\caption{\small Numerical comparison of the new bounds 
\eqref{blmeamajorbl} and \eqref{blmervimajorbl}
with the relevant bounds
\eqref{blmeamajorblui} and \eqref{blmervimajorblui}
based on \cite{lz2015} accompanying
the block Chebyshev and Lanczos methods;
see Example 1 with $p=3$, $\tau=\{1,2,3\}$, $t=3$ and $i=3$.}
\label{fig1}
\end{figure}

Similarly to \cite[Subsection 4.2]{s1980} and \cite[Example 7.3]{lz2015}, we consider
the diagonal matrix $A=\mbox{diag}(\la_1,\ldots,\la_n)$ with $n=900$ and the eigenvalues
\[\la_1=2, \quad \la_2=1.6, \quad \la_3=1.4, \quad 
 \la_j=1-(j-3)/n \ \ \mbox{for} \ \ j=4,\ldots,n.\]
Combining this with the settings from Subsection 2.1 and Lemma \ref{lm:auxvec},
the invariant subspace $\mathcal{X}$ associated with the $p$ largest eigenvalues
is spanned by the first $p$ columns of the identity matrix $I_n$. We set $p=3$ and construct
the initial subspace $\mathcal{Y}\in\R^{n \times p}$ of the block Lanczos method
using a random matrix \,\verb|Y=[orth(randn(p,p)); randn(n-p,p)]|\,
such that the principal angles between $\mathcal{X}$ and $\mathcal{Y}$
are evenly distributed. We note that the initial subspace in the related examples
in \cite{s1980,lz2015} is spanned by a fixed matrix
\,$Y=[V;\,\ldots;\,V]$\, with \,$V=[1,\ 1,\ 1;\ 1,\ 0,\ -2;\ 1,\ {-}1,\ 1]$\,
for which the principal angles between $\mathcal{X}$ and $\mathcal{Y}$ are actually equal.
We use $1000$ randomly constructed initial subspaces instead.
For each of them, the corresponding block Krylov subspace
$\K=\Y+A\Y+\cdots+A^{k-1}\Y$ is determined up to $k=15$
where full orthogonalization is used in order to reduce instability.
Indeed, the target eigenvalues $\la_1,\la_2,\la_3$ in this classical example
are well separated, and the gap $\la_3-\la_4\approx0.4$ is sufficiently large
to ensure meaningful Chebyshev terms in the bounds. Figure \ref{fig1} shows
the numerical comparison between our new bounds
\eqref{blmeamajorbl}, \eqref{blmervimajorbl}
and the reformulated bounds \eqref{blmeamajorblui}, \eqref{blmervimajorblui}
based on \cite{lz2015}. In the left subfigure, the multi-angle bounds \eqref{blmeamajorbl}
and \eqref{blmeamajorblui} are compared in the special form
\[\textstyle\sum_{j=1}^t\tan\theta_j(\X_{\tau},\K)
 \le\sum_{j=1}^t\sigma\tan\theta_j(\X_{\tau},\Y_{\tau})
 \quad\mbox{with}\quad\sigma\in[\sigma_{i_t},\ldots,\sigma_{i_1}]
 \quad\mbox{or}\quad\sigma=\sigma_{i_t}.\]
The ``Lanczos'' curve illustrates the mean value (cf.~\cite{lz2015})
of \,$\textstyle\sum_{j=1}^t\tan\theta_j(\X_{\tau},\K)$\,
among $1000$ samples for $\tau=\{1,2,3\}$ and $t=3$ for each $k\in\{1,\ldots,15\}$.
The ``Chebyshev'' curve presents the corresponding data
determined for $f(A)\Y$ instead of $\K$ by using
the shifted Chebyshev polynomial defined in \eqref{scp}.
The other curves contain the mean values of bounds.
The accuracy of \eqref{blmeamajorblui} is at the same level as that observed
in \cite[Table 4]{lz2015} concerning the measure
$\big(\sum_{j=1}^t\sin^2\theta_j(\X_{\tau},\K)\big)^{1/2}$.
The improvement achieved by \eqref{blmeamajorbl} is evident, and is
also observed in a comparison with respect to maxima instead of mean values.
In the right subfigure, we evaluate
\[\textstyle\sum_{j=1}^i\dfrac{\la_j-\psi_j}{\la_1-\la_n}
 \le\sum_{j=1}^i\sigma^2\tan^2\theta_j(\X_i,\Y_i)
 \quad\mbox{with}\quad\sigma\in[\sigma_i,\ldots,\sigma_1]
 \quad\mbox{or}\quad\sigma=\sigma_i\]
to compare the Ritz value bounds \eqref{blmervimajorbl}
and \eqref{blmervimajorblui} for $i=3$. Therein the denominator $\psi_j-\la_n$
in \eqref{blmervimajorbl} is simplified as $\la_1-\la_n$.
This reduces the accuracy of \eqref{blmervimajorbl}, but only slightly
after several iteration steps. A possible overestimation
in the first iteration steps is avoided by additionally using the trivial bound
$\sum_{j=1}^i\frac{\la_j-\psi_j}{\la_1-\la_n} \le i$.
We observe that \eqref{blmervimajorbl} is more accurate
than \eqref{blmervimajorblui}. However, their relative accuracy is slightly
worse than that of \eqref{blmeamajorbl} and \eqref{blmeamajorblui}.
This reflects the fact that the Ritz value bounds are essentially derived by
combining the multi-angle bounds with further inequalities
which are not necessarily sharp at the same time.
Furthermore, the reader may refer to \cite[Section 7]{lz2015}
and \cite[Section 5]{z2018} that discuss the drawbacks
of the classical bounds from \cite{s1980}.

\subsection{Example 2}
\begin{figure}[htbp]
\begin{center}
\includegraphics[width=0.9\textwidth]{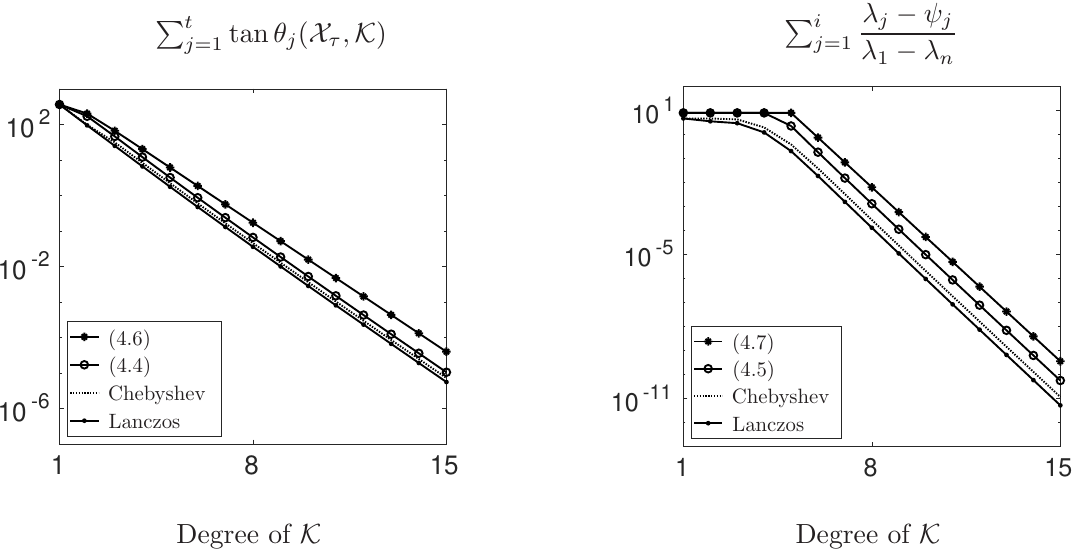}\quad
\end{center}
\par\vskip -2ex
\caption{\small Numerical comparison of the new bounds 
\eqref{blmeamajorbl} and \eqref{blmervimajorbl}
with the relevant bounds
\eqref{blmeamajorblui} and \eqref{blmervimajorblui}
based on \cite{lz2015} accompanying
the block Chebyshev and Lanczos methods;
see Example 2 with $p=9$, $\tau=\{3,\ldots,8\}$, $t=6$ and $i=8$.}
\label{fig2}
\end{figure}

We use the diagonal matrix $A=\mbox{diag}(\la_1,\ldots,\la_n)$
with $n=3600$ and the eigenvalues
\[\begin{split}
 &\la_1=2.05, \quad \la_2=2, \quad \la_3=1.95, \quad
   \la_4=1.65, \quad \la_5=1.6, \quad \la_6=1.55,\\
 &\la_7=1.45, \quad \la_8=1.4, \quad \la_9=1.35, \quad
   \la_j=1-(j-9)/n \ \ \mbox{for} \ \ j=10,\ldots,n.
\end{split}\]
The experiment introduced in Example 1 is run for this matrix together with
$p=9$, $\tau=\{3,\ldots,8\}$, $t=6$ and $i=8$; see Figure \ref{fig2}.
The smaller distances between target eigenvalues do not deteriorate
the convergence rates due to the large gap $\la_9-\la_{10}\approx0.35$.
This also results in suitable Chebyshev terms
so that the slopes of the bound curves accurately reflect
the convergence rates and the cluster robustness, at least in the final phase.
Our new bounds \eqref{blmeamajorbl} and \eqref{blmervimajorbl}
provide visible improvements
in comparison to \eqref{blmeamajorblui} and \eqref{blmervimajorblui}
based on \cite{lz2015}. The corresponding bounds
from \cite{lz2015} already dramatically improve the classical bounds from \cite{s1980}
involving the ``bulky'' factors; see \cite[Example 7.2]{lz2015}.

\section*{Conclusions}

Majorization techniques exploring relations between Ritz value errors
and principal angles are extended to block signal filters
and subspace iterates of block eigensolvers.
The majorization-type analysis of tangents of principal angles from \cite{ka10}
is improved by using auxiliary vectors from the classical analysis
of the block power method by Rutishauser \cite{r1969}. This leads to
novel majorization bounds for the Rayleigh-Ritz method applied to
the final iterative subspace in terms of principal angles between
an initial subspace and a target invariant subspace. Our results improve
the existing approaches from \cite{s1980,lz2015}
that are used for the block Lanczos method.
Our majorization technique is especially advantageous in situations
where the concerned principal angles are evenly distributed,
which is common for random initial subspaces. 
Our forthcoming work will
focus on angle-free bounds for Ritz value errors which are applicable to
restarted iterations, such as the stochastic block descent.
Simultaneous consideration
of several Ritz values approximating clustered eigenvalues enables
meaningful bounds and has potential for investigating
block preconditioned eigensolvers.  A future aim is to cover the popular 
locally optimal block preconditioned conjugate gradient (LOBPCG) method
\cite{k01,klao07}, taking advantage of techniques developed in \cite{o2006,zn2019}.

\section{Appendix: Detailed proofs}

\subsection{Proof of Lemma \ref{lm:tangentform}}\label{p:tangentform}

By using the orthonormal basis matrix $U=\widetilde{U}G$
with $G=(\widetilde{U}^H\widetilde{U})^{-1/2}\in\C^{s \times s}$,
the cosine values of the principal angles 
$\theta_1\ge\cdots\ge\theta_s$ from ${\U}$ to ${\V}$
are given by the singular values of $V^HU\in\C^{t \times s}$.
More precisely, a standard singular value decomposition
$V^HU=W \Sigma Z^H$
provides unitary matrices $W\in\C^{t \times t}$ and $Z\in\C^{s \times s}$
together with the diagonal (rectangle) matrix
\[\Sigma=\begin{bmatrix}D\\O\end{bmatrix}\in\R^{t \times s}
 \quad\mbox{consisting of}\quad
 D=\diag\big(\cos(\theta_s),\ldots,\cos(\theta_1)\big)\in\R^{s \times s}\]
and the zero matrix $O\in\R^{(t-s) \times s}$.
Then $\widehat{V}=VW$ and $\widehat{U}=UZ=\widetilde{U}GZ$
are orthonormal basis matrices of $\U$ and $\V$
with the property $\Sigma=\widehat{V}^H\widehat{U}=W^HV^H\widetilde{U}GZ$.
Since $V^H\widetilde{U}$ has full rank, so does $\Sigma$.
In addition, it holds that
\[\begin{split}
 \QVP^H\widehat{U}\Sigma^{\dag}
 &=\QVP^H\widehat{U}(\widehat{V}^H\widehat{U})^{\dag}
 =\QVP^H\widetilde{U}GZ(W^HV^H\widetilde{U}GZ)^{\dag}
 =\QVP^H\widetilde{U}GZ(GZ)^{\dag}(V^H\widetilde{U})^{\dag}(W^H)^{\dag}\\[1ex]
 &=\QVP^H\widetilde{U}GZ(GZ)^{-1}(V^H\widetilde{U})^{\dag}(W^H)^{-1}
 =\big(\QVP^H\widetilde{U}(V^H\widetilde{U})^{\dag}\big)W
\end{split}\]
so that the singular values of $\QVP^H\widehat{U}\Sigma^{\dag}$
coincide with those of $\QVP^H\widetilde{U}(V^H\widetilde{U})^{\dag}$.
Their squared values are the $\min\{t,\,n{-}t\}$
largest eigenvalues of the $t{\times}t$ matrix
\begin{equation}\label{svd}
\begin{split}
 (\QVP^H\widehat{U}\Sigma^{\dag})^H\QVP^H\widehat{U}\Sigma^{\dag}
 &=(\Sigma^{\dag})^H\widehat{U}^H\QVP\QVP^H\widehat{U}\Sigma^{\dag}
 =(\Sigma^{\dag})^H\widehat{U}^H(I_n-\widehat{V}\widehat{V}^H)
 \widehat{U}\Sigma^{\dag}\\[1ex]
 &=(\Sigma^{\dag})^H\widehat{U}^H\widehat{U}\Sigma^{\dag}
 -(\Sigma^{\dag})^H\widehat{U}^H\widehat{V}\widehat{V}^H
 \widehat{U}\Sigma^{\dag}\\[1ex]
 &=(\Sigma^{\dag})^H\Sigma^{\dag}-(\Sigma^{\dag})^H\Sigma^H\Sigma\Sigma^{\dag}
 =(\Sigma^{\dag})^H\Sigma^{\dag}-(\Sigma\Sigma^{\dag})^H(\Sigma\Sigma^{\dag}).
\end{split}
\end{equation}
Therein $\QVP\QVP^H$ is the orthogonal projector on $\VP$ and thus
coincides with $I_n-\widehat{V}\widehat{V}^H$
for the orthogonal projector $\widehat{V}\widehat{V}^H$ on $\V$.
Moreover, $\Sigma^{\dag}=[D^{-1}\ O^H]$ holds as $\Sigma$ has full rank.
Then \eqref{svd} implies
\[(\QVP^H\widehat{U}\Sigma^{\dag})^H\QVP^H\widehat{U}\Sigma^{\dag}
 =\begin{bmatrix}D^{-2}&\\&\widetilde{O}\end{bmatrix}
 -\begin{bmatrix}I_s&\\&\widetilde{O}\end{bmatrix}^H
 \begin{bmatrix}I_s&\\&\widetilde{O}\end{bmatrix}
 =\begin{bmatrix}D^{-2}-I_s&\\&\widetilde{O}\end{bmatrix}\]
with the zero matrix $\widetilde{O}=OO^H$. Thus the $s$ largest singular values of
$\QVP^H\widetilde{U}(V^H\widetilde{U})^{\dag}$ or $\QVP^H\widehat{U}\Sigma^{\dag}$
coincide with the square roots of the diagonal entries of the diagonal matrix
\[D^{-2}-I_s=\diag\Big(\big(\cos(\theta_s)\big)^{-2}-1,\,\ldots,\,
 \big(\cos(\theta_1)\big)^{-2}-1\Big)
 =\diag\big(\tan^2(\theta_s),\ldots,\tan^2(\theta_1)\big),\]
i.e., the tangent values of the principal angles from ${\U}$ to ${\V}$.
The largest singular value also coincides with the $2$-norm so that
$\|\QVP^H\widetilde{U}(V^H\widetilde{U})^{\dag}\|
=\tan(\theta_1)=\tan\angle(\U,\V)$.
\hfill$\Box$

\subsection{Proof of Lemma \ref{lm:blme}}\label{p:blme}

Applying Lemma \ref{lm:tangentform} to
\[V=x_i,\quad\QVP=[x_1,\ldots,x_{i-1},x_{i+1},\ldots,x_n],\quad\widetilde{U}=y_i\]
shows that the only singular value of the $(n{-}1){\times}1$ matrix
\[w=[x_1,\ldots,x_{i-1},x_{i+1},\ldots,x_n]^Hy_i(x_i^Hy_i)^{\dag}
 \stackrel{\eqref{biorth}}{=}[0;\,\ldots;\,0;\,x_{p+1}^Hy_i;\,\ldots;\,x_n^Hy_i].\]
coincides with the tangent value of the only principal angle
from $\span\{y_i\}$ to $\span\{x_i\}$. Then
\[\|w\|=|\tan\angle(y_i,x_i)|=|\tan\angle(x_i,y_i)|=\tan\angle(x_i,y_i)\]
holds since $x_i^Hy_i=1>0$ leads to $\tan\angle(x_i,y_i)>0$.
In addition, the common assumption \eqref{ca} ensures $f(\la_i)\neq0$
and $x_i^Hf(A)y_i=\big(\big(f(A)\big)^Hx_i\big)^Hy_i=\big(\overline{f(\la_i)}x_i\big)^Hy_i
=f(\la_i)x_i^Hy_i=f(\la_i)\neq0$.
Thus $f(A)y_i$ is a nonzero vector. Applying Lemma \ref{lm:tangentform} to
\[V=x_i,\quad\QVP=[x_1,\ldots,x_{i-1},x_{i+1},\ldots,x_n],\quad\widetilde{U}=f(A)y_i\]
analogously implies $\big|\tan\angle\big(x_i,f(A)y_i\big)\big|=\|w'\|$ for
\[\begin{split}
 w' \,\ = \,\ &[x_1,\ldots,x_{i-1},x_{i+1},\ldots,x_n]^Hf(A)y_i\big(x_i^Hf(A)y_i\big)^{\dag}\\[1ex]
 \stackrel{\eqref{biorth}}{=}
 &\big[0;\,\ldots;\,0;\,f(\la_{p+1})x_{p+1}^Hy_i;\,\ldots;\,f(\la_n)x_n^Hy_i\big]
 \big(f(\la_i)\big)^{-1}.
\end{split}\]
Then a simple comparison between $\|w'\|$ and $\|w\|$ shows
\[\big|\tan\angle\big(x_i,f(A)y_i\big)\big|
 \le\frac{\max_{j\in\{p+1,\ldots,n\}}|f(\la_j)|}{|f(\la_i)|}\tan\angle(x_i,y_i).\]
This implies the first inequality in \eqref{blmea} according to
$f(A)y_i\in f(A)\Y=\Y'$.

Moreover, the term $\tan\angle(x_i,y_i)$
can be enlarged by $\tan\angle(\X,\Y)$
for eliminating the auxiliary vector $y_i$; cf.~\eqref{blme}.
This modification is enabled by projection arguments
or by applying Lemma \ref{lm:tangentform} to
\[V=X=[x_1,\ldots,x_p],\quad\QVP=[x_{p+1},\ldots,x_n],\quad\widetilde{U}=y_i.\]
Therein $\tan\angle(y_i,\X)=\|W\|$ holds for
\[W=[x_{p+1},\ldots,x_n]^Hy_i\big([x_1,\ldots,x_p]^Hy_i\big)^{\dag}
 =\widetilde{w}e_i^{\dag}=\widetilde{w}e_i^H\]
with $\widetilde{w}=[x_{p+1}^Hy_i;\,\ldots;\,x_n^Hy_i]$ and 
the $i$th column $e_i$ of the identity matrix $I_p$. Then
\[\|W\|=\|\widetilde{w}e_i^H\|=\|\widetilde{w}\|\|e_i\|
 =\|\widetilde{w}\|=\|w\|\]
holds by comparing the components of $w$ and $\widetilde{w}$. Thus
\[\tan\angle(x_i,y_i)=\|w\|=\|W\|=\tan\angle(y_i,\X)
 \le\tan\angle(\Y,\X)=\tan\angle(\X,\Y)\]
leads to the second inequality in \eqref{blmea}.
\hfill$\Box$

Bound \eqref{blmea} can also be shown in an elementary way;
cf.~\cite[Theorem 1]{z2018} concerning the specified function \eqref{fsp}.
The proof of Lemma \ref{lm:blme} aims at
motivating that of Theorem \ref{thm:blmemajor}
where majorization bounds are to be derived.

\subsection{Proof of Lemma \ref{lm:blmervi}}\label{p:blmervi}

The common assumption \eqref{ca} ensures $f(\la_j)\neq0$
for each $j\in\{1,\ldots,i\}$. The $i{\times}i$ matrix
\[X_i^HY'_i=X_i^Hf(A)Y_i=\diag\big(f(\la_1),\ldots,f(\la_i)\big)X_i^HY_i
 \stackrel{\eqref{biorth}}{=}\diag\big(f(\la_1),\ldots,f(\la_i)\big)\]
is then invertible. Thus $Y'_i$ has rank $i$, and $\dim\Y'_i=i$.
Moreover, the relation $\Y'=f(A)\,\span\{Y\}=\span\{f(A)Y\}
=\span\{f(A)Y_p\}=\span\{Y'_p\}=\Y'_p$ holds
so that $\dim\Y'=p$.

In order to show \eqref{blmervia},
the property $\dim\Y'_i=i=\dim\X_i$ enables the representation
\[\tan\angle(\X_i,\Y'_i)=\tan\angle(\Y'_i,\X_i)
 =\max_{y'\in\Y'_i{\setminus}\{0\}}\tan\angle(y',\X_i)
 =\tan\angle(\widehat{y}',\X_i)\]
with a maximizer $\widehat{y}'\in\Y'_i{\setminus}\{0\}$.
Because $\Y'_i=\span\{f(A)Y_i\}=f(A)\Y_i$,
we can represent $\widehat{y}'$ by $\widehat{y}'=f(A)\widehat{y}$
with a certain $\widehat{y}\in\Y_i{\setminus}\{0\}$.
Applying Lemma \ref{lm:tangentform} to
\[V=X_i,\quad\QVP=[x_{i+1},\ldots,x_n],\quad\widetilde{U}=\widehat{y}'=f(A)\widehat{y}\]
(there is only one principal angle in this case)
yields $\tan\angle(\widehat{y}',\X_i)=\|W'\|$ for
\[\begin{split}
 & W'=[x_{i+1},\ldots,x_n]^Hf(A)\widehat{y}\,\big(X_i^Hf(A)\widehat{y}\big)^{\dag}
  =w'{e'}^{\dag}\\[1ex]
 &\mbox{with} \quad w'=\big[0;\,\ldots;\,0;\,f(\la_{p+1})x_{p+1}^H\widehat{y};
 \,\ldots;\,f(\la_n)x_n^H\widehat{y}\big]\\[1ex]
 &\mbox{and} \quad \
 e'=\big[f(\la_1)x_1^H\widehat{y};\,\ldots;\,f(\la_i)x_i^H\widehat{y}\big].
\end{split}\]
The vector $w'$ initially consists of \,$x_j^Hf(A)\widehat{y}=f(\la_j)x_j^H\widehat{y}$\,
for $j=i{+}1,\ldots,n$, but the first components up to the index $p$ are simply zero
due to \eqref{biorth} and that $\widehat{y}$ belongs to $\Y_i=\span\{y_1,\ldots,y_i\}$. 
Moreover, $\widehat{y}$ can be represented by
$\widehat{y}=Y_ig$ with a certain $g\in\C^i{\setminus}\{0\}$ so that
\[e'=\diag\big(f(\la_1),\ldots,f(\la_i)\big)X_i^HY_ig
 \stackrel{\eqref{biorth}}{=}\diag\big(f(\la_1),\ldots,f(\la_i)\big)g\neq0.\]
Then ${e'}^{\dag}={e'}^H/\|e'\|^2$, and
\[\tan\angle(\widehat{y}',\X_i)
 =\|W'\|=\big\|w'{e'}^H/\|e'\|^2\big\|=\frac{\|w'{e'}^H\|}{\|e'\|^2}
 =\frac{\|w'\|\|e'\|}{\|e'\|^2}=\frac{\|w'\|}{\|e'\|}.\]
Analogously, $\tan\angle(\widehat{y},\X_i)$ has the representation
\[\begin{split}\tan\angle(\widehat{y},\X_i)=\frac{\|w\|}{\|e\|}
 \quad & \mbox{with} \quad
 w=[0;\,\ldots;\,0;\,x_{p+1}^H\widehat{y};\,\ldots;\,x_n^H\widehat{y}]\\
 & \mbox{and} \quad \
 e=[x_1^H\widehat{y};\,\ldots;\,x_i^H\widehat{y}].
\end{split}\]
Summarizing the above together with
\[\|w'\|\le\big(\max_{j\in\{p+1,\ldots,n\}}|f(\la_j)|\big)\|w\|,\quad
 \|e'\|\ge\big(\min_{j\in\{1,\ldots,i\}}|f(\la_j)|\big)\|e\|\]
and $\tan\angle(\widehat{y},\X_i)\le\tan\angle(\Y_i,\X_i)=\tan\angle(\X_i,\Y_i)$
yields \eqref{blmervia}.
\hfill$\Box$

The above proof focuses on the rank-$1$ matrices $w'{e'}^{\dag}$ and $we^{\dag}$.
An alternative proof using rank-$i$ matrices is involved in the corresponding
majorization-type analysis; see the proofs of Theorem \ref{thm:blmemajor}
and Lemma \ref{lm:blmervimajor}.

\subsection{Proof of Lemma \ref{lm:blmervi1}}\label{p:blmervi1}

We use the smallest Ritz value of $A$ in $\Y'_i$ and an associated
normalized Ritz vector which are denoted by $\eta'$ and $y'$. Then
\[\Y'_i\subseteq\Y'\subset\C^n\quad\Rightarrow\quad
 \la_i\ge\eta'_i\ge\eta'\ge\la_n\]
according to the Courant-Fischer principles and $\dim\Y'_i=i$.
Moreover, \eqref{blmervia1} can be implied by
\begin{equation}\label{blmervia1a}
 \frac{\la_i-\eta'_i}{\eta'_i-\la_n}\le\tan^2\angle(y',\X_i)
\end{equation}
since $y'$ is a nonzero vector in $\Y'_i$ and thus
$\tan^2\angle(y',\X_i)\le\tan^2\angle(\Y'_i,\X_i)=\tan^2\angle(\X_i,\Y'_i)$.
An elementary trigonometric transformation shows that
\eqref{blmervia1a} is equivalent to
\begin{equation}\label{blmervia2}
 \cos^2\angle(y',\X_i)\le\frac{\eta'_i-\la_n}{\la_i-\la_n}.
\end{equation}
Thus we only need to show \eqref{blmervia2}.
Therein we use the standard definition of principal angles
which indicates that $\cos\angle(y',\X_i)$ coincides with the only singular value of
$X_i^Hy'$, i.e., $\|X_i^Hy'\|$. In addition, we consider the Hermitian matrices
\[\widetilde{A}=A-\la_nI_n\quad\mbox{and}\quad
 \widetilde{D}=\diag(\la_1-\la_n,\,\ldots,\,\la_i-\la_n).\]
Evidently, $\widetilde{A}$ is positive semidefinite,
and $\widetilde{D}$ is positive definite because
$\la_1\ge\cdots\ge\la_i\ge\la_p>\la_{p+1}\ge\la_n$. By using their square roots
$\widetilde{A}^{1/2}$ and $\widetilde{D}^{1/2}$, we get
$\widetilde{A}^{1/2}X_i=X_i\widetilde{D}^{1/2}$ so that
\[\begin{split}
 &X_i^Hy'=\widetilde{D}^{-1/2}\widetilde{D}^{1/2}X_i^Hy'
 =\widetilde{D}^{-1/2}\big(X_i\widetilde{D}^{1/2}\big)^Hy'
 =\widetilde{D}^{-1/2}\big(\widetilde{A}^{1/2}X_i\big)^Hy'
 =\widetilde{D}^{-1/2}X_i^H\widetilde{A}^{1/2}y'\\[1ex]
 \Rightarrow\quad&\cos\angle(y',\X_i)
 =\|X_i^Hy'\|=\|\widetilde{D}^{-1/2}X_i^H\widetilde{A}^{1/2}y'\|
 \le\|\widetilde{D}^{-1/2}\|\|X_i^H\widetilde{A}^{1/2}y'\|.
\end{split}\]
Combining this with
\[\begin{split}
 \|\widetilde{D}^{-1/2}\|
 &=\max_{j=1,\ldots,i}\big|(\la_j-\la_n)^{-1/2}\big|=(\la_i-\la_n)^{-1/2}
 \quad\mbox{and}\\[1ex]
 \|X_i^H\widetilde{A}^{1/2}y'\|
 &\le\|X_i^H\|\|\widetilde{A}^{1/2}y'\|
 =\|\widetilde{A}^{1/2}y'\|=\big({y'}^H\widetilde{A}y'\big)^{1/2}\\[1ex]
 &=\big({y'}^HAy'-\la_n\big)^{1/2}
 =(\eta'-\la_n)^{1/2}\le(\eta'_i-\la_n)^{1/2}
\end{split}\]
yields \eqref{blmervia2}.
\hfill$\Box$

The proof of Lemma \ref{lm:blmervi1} uses inequalities
with respect to the $2$-norm, i.e., the largest singular value.
In the corresponding majorization-type analysis,
we consider tuples of singular values;
see the proof of Lemma \ref{lm:blmervimajor1}.

\subsection{Extending \eqref{blmerviabi1} as \eqref{blmerviabi}}\label{p:ext}

For this extension, we apply Lemma \ref{lm:tangentform} to
\[V=X_i=[x_1,\ldots,x_i],\quad\QVP=[x_{i+1},\ldots,x_n],
 \quad\widetilde{U}=Y_i=[y_1,\ldots,y_i].\]
This leads to $\tan\angle(\X_i,\Y_i)=\tan\angle(\Y_i,\X_i)=\|W\|$ for
\[W=[x_{i+1},\ldots,x_n]^HY_i\big(X_i^HY_i\big)^{\dag}
 \stackrel{\eqref{biorth}}{=}[O;\,x_{p+1}^HY_i;\,\ldots;\,x_n^HY_i]\,I_i^{\dag}
 =\begin{bmatrix}O\\G\end{bmatrix}\]
with $G=[x_{p+1}^HY_i;\,\ldots;\,x_n^HY_i]\in\C^{(n-p) \times i}$
and the zero matrix $O\in\R^{(p-i) \times i}$.
Moreover, applying Lemma \ref{lm:tangentform} to
\[V=X=[x_1,\ldots,x_p],\quad\QVP=[x_{p+1},\ldots,x_n],
 \quad\widetilde{U}=Y=[y_1,\ldots,y_p]\]
yields $\tan\angle(\X,\Y)=\tan\angle(\Y,\X)=\|\widehat{G}\|$ for
\[\widehat{G}=[x_{p+1},\ldots,x_n]^HY\big(X^HY\big)^{\dag}
 \stackrel{\eqref{biorth}}{=}[x_{p+1}^HY;\,\ldots;\,x_n^HY]\,I_p^{\dag}
 =[x_{p+1}^HY;\,\ldots;\,x_n^HY]\in\C^{(n-p) \times p}.\]
In summary, it holds that
\[\tan\angle(\X_i,\Y_i)=\|W\|=\|G\|
 =\|\widehat{G}E_i\|\le\|\widehat{G}\|=\tan\angle(\X,\Y)\]
where $E_i\in\R^{p\times i}$ consists of the first $i$ columns of $I_p$.

\subsection{Proof of Theorem \ref{thm:blmemajor}}\label{p:blmemajor}

Applying Lemma \ref{lm:tangentform} to
\[V=X_{\tau}=[x_{i_1},\ldots,x_{i_t}],\quad
 \QVP=[x_j\ \mbox{for} \ j\in\{1,\ldots,n\}{\setminus}\tau],\quad
 \widetilde{U}=Y_{\tau}=[y_{i_1},\ldots,y_{i_t}]\]
and $s=t$ generates the $(n{-}t){\times}t$ matrix
\[W=[x_j\ \mbox{for} \ j\in\{1,\ldots,n\}{\setminus}\tau]^H
 Y_{\tau}(X_{\tau}^HY_{\tau})^{\dag}
 \stackrel{\eqref{biorth}}{=}[O;
 \,x_{p+1}^HY_{\tau};\,\ldots;\,x_n^HY_{\tau}]\]
with the zero matrix $O\in\R^{(p-t) \times t}$
where $(X_{\tau}^HY_{\tau})^{\dag}=I_t^{\dag}=I_t$ is eliminated.
Moreover, $t \le p \ll n$ holds according to Subsection 2.1
and ensures $t \le n{-}t$ so that $W$ has $t$ singular values.
By using the notations from Subsection 2.3, the statement
in Lemma \ref{lm:tangentform} yields
\[S(W)=\tan\Theta(\Y_{\tau},\X_{\tau})=\tan\Theta(\X_{\tau},\Y_{\tau}).\]
In addition, the common assumption \eqref{ca} ensures
$f(\la_j)\neq0 \ \ \forall\ j\in\tau$ so that
\[D_\tau{\,=\,}X_{\tau}^Hf(A)Y_{\tau}
 =\diag\big(f(\la_{i_1}),\ldots,f(\la_{i_t})\big)X_{\tau}^HY_{\tau}
 \stackrel{\eqref{biorth}}{=}
 \diag\big(f(\la_{i_1}),\ldots,f(\la_{i_t})\big)\]
is an invertible diagonal matrix.
Consequently, the $n{\times}t$ matrix
$f(A)Y_{\tau}$ has rank $t$, and the subspace
$\Y'_{\tau}=\span\{f(A)Y_{\tau}\}$ has dimension $t$.
Applying Lemma \ref{lm:tangentform} to
\[V=X_{\tau},\quad
 \QVP=[x_j\ \mbox{for} \ j\in\{1,\ldots,n\}{\setminus}\tau],\quad
 \widetilde{U}=f(A)Y_{\tau}\]
implies \,$S(W')=\tan\Theta(\X_{\tau},\Y'_{\tau})$\,
for the $(n{-}t){\times}t$ matrix
\[\begin{split}
 W' \,\ = \,\ &[x_j\ \mbox{for} \ j\in\{1,\ldots,n\}{\setminus}\tau]^Hf(A)Y_{\tau}
 \big(X_{\tau}^Hf(A)Y_{\tau}\big)^{\dag}\\[1ex]
 \stackrel{\eqref{biorth}}{=}
 &\big[O;\,f(\la_{p+1})x_{p+1}^HY_{\tau};\,\ldots;\,
 f(\la_n)x_n^HY_{\tau}\big]D_\tau^{\dag}\\[1ex]
 = \,\ &\widehat{D}\,WD_\tau^{-1}\quad\mbox{with}\quad
 \widehat{D}=\diag\big(0,\ldots,0,f(\la_{p+1}),\ldots,f(\la_n)\big).
\end{split}\]
Furthermore, applying Lemma \ref{lm:major} to
\[B_1=\widehat{D},\quad B_2=W,\quad B_3=D_\tau^{-1}\]
and $c=1$ leads to
\[\tan\Theta(\X_{\tau},\Y'_{\tau})=S(W')
 =S_t(W')\prec_w S_t(\widehat{D})S_t(W)S_t(D_\tau^{-1})
 =S(D_\tau^{-1})S_t(\widehat{D})S(W).\]
Combining this with
\[\begin{split}
 &S(D_\tau^{-1})=\left[\big|\big(f(\la_{i_1})\big)^{-1}\big|,\ldots,
 \big|\big(f(\la_{i_t})\big)^{-1}\big|\right]^{\downarrow}=\Phi_{\tau},\\
 &S(\widehat{D})=\big[0,\ldots,0,
 |f(\la_{p+1})|,\ldots,|f(\la_n)|\big]^{\downarrow}\quad\Rightarrow\quad
 S_t(\widehat{D})=\widehat{\Phi}_t
\end{split}\]
and $S(W)=\tan\Theta(\X_{\tau},\Y_{\tau})$ yields 
\begin{equation}\label{blmeamajor1}
 \tan\Theta(\X_{\tau},\Y'_{\tau}) \prec_w
 \Phi_{\tau}\,\widehat{\Phi}_t\,
 \tan\Theta(\X_{\tau},\Y_{\tau}).
\end{equation}
Then \eqref{blmeamajor} is obtained by using the tuple inequality
$\tan\Theta(\X_{\tau},\Y')\le\tan\Theta(\X_{\tau},\Y'_{\tau})$
based on $\Y'_{\tau}\subseteq\span\{f(A)Y\}=\Y'$. This inequality
can be shown in the equivalent form
$\big(\cos^2\Theta(\X_{\tau},\Y')\big)^{\downarrow}
 \ge\big(\cos^2\Theta(\X_{\tau},\Y'_{\tau})\big)^{\downarrow}$
as follows. An arbitrary orthonormal basis matrix $\widetilde{Y}$
of $\Y'_{\tau}$ can be extended as an orthonormal basis matrix
$[\widetilde{Y} \,\ \widehat{Y}]$ of $\Y'$ so that
\[\begin{split}
 \big(\cos^2\Theta(\X_{\tau},\Y')\big)^{\downarrow}
 &=S^2\big([\widetilde{Y} \,\ \widehat{Y}]^HX_{\tau}\big)
 =\Lambda\big(\big([\widetilde{Y} \,\ \widehat{Y}]^HX_{\tau}\big)^H
 [\widetilde{Y} \,\ \widehat{Y}]^HX_{\tau}\big)\\[1ex]
 &=\Lambda\big(X_{\tau}^H\widetilde{Y}\widetilde{Y}^HX_{\tau}
 +X_{\tau}^H\widehat{Y}\widehat{Y}^HX_{\tau}\big)
 \ge\Lambda\big(X_{\tau}^H\widetilde{Y}\widetilde{Y}^HX_{\tau}\big)
 =\big(\cos^2\Theta(\X_{\tau},\Y'_{\tau})\big)^{\downarrow}.
\end{split}\]
Therein the intermediate inequality uses the Weyl's inequality
and the fact that $X_{\tau}^H\widehat{Y}\widehat{Y}^HX_{\tau}$
is Hermitian positive semidefinite.
\hfill$\Box$

\subsection{Derivation of \eqref{blmeamajorab}}\label{p:blmeamajorab}

For deriving \eqref{blmeamajorab},
we use the given basis matrix $\widetilde{U}$ of $\U$ together with
the orthogonal projectors $VV^H$ on $\V$ so that
\[V^HF\widetilde{U}=(F^HV)^H\widetilde{U}
 =\big((VV^H)F^HV\big)^H\widetilde{U}=(V^HFV)V^H\widetilde{U}.\]
Analogously, $\QVP^HF\widetilde{U}=(\QVP^HF\QVP)\QVP^H\widetilde{U}$.
Moreover, $V^HFV\in\C^{t \times t}$ is invertible
and $V^H\widetilde{U}\in\C^{t \times s}$
has full rank by setting, thus $V^HF\widetilde{U}$ has full rank.
Then $F\widetilde{U}$ also has full rank
and can be treated as a basis matrix of $F\U$ with Lemma \ref{lm:tangentform}.
By using the matrices
\[G=V^HFV, \quad \widetilde{W}=V^H\widetilde{U}, \quad
 D=\big(\widetilde{W}^H\widetilde{W}\big)^{1/2}
 \quad\mbox{and}\quad W=\widetilde{W}D^{-1},\]
it holds that
\[(V^HF\widetilde{U})^{\dag}=(G\,\widetilde{W})^{\dag}
 =\widetilde{W}^{\dag}\widetilde{W}(G\,WD)^{\dag}
 =\widetilde{W}^{\dag}\widetilde{W}D^{-1}(G\,W)^{\dag}
 =(V^H\widetilde{U})^{\dag}W(G\,W)^{\dag}.\]
In summary, we get
\[\QVP^HF\widetilde{U}(V^HF\widetilde{U})^{\dag}
 =(\QVP^HF\QVP)\big(\QVP^H\widetilde{U}(V^H\widetilde{U})^{\dag}\big)
 \big(W(G\,W)^{\dag}\big).\]
Subsequently, we consider
the singular values $\alpha_1\le\cdots\le\alpha_t$ of $G$
and the singular values $\beta_1\le\cdots\le\beta_s$ of $G\,W$.
These values are evidently all positive, and it holds that
$\alpha_i\le\beta_i$, $i=1,\ldots,s$\, according to 
$W^HW=I_s$ and the Courant-Fischer principles.
This leads to the tuple inequality
\[S_s\big(W(G\,W)^{\dag}\big)=S\big((G\,W)^{\dag}\big)
 =[\beta_1^{-1},\ldots,\beta_s^{-1}]\le[\alpha_1^{-1},\ldots,\alpha_s^{-1}]
 =S_s\big(G^{-1}\big)=S_s\big((V^HFV)^{-1}\big)\]
which can easily be combined with Lemma \ref{lm:major} for showing
\[S_s\big(\QVP^HF\widetilde{U}(V^HF\widetilde{U})^{\dag}\big)
 \prec_w S_s\big(\QVP^HF\QVP\big)\,
 S_s\big(\QVP^H\widetilde{U}(V^H\widetilde{U})^{\dag}\big)\,
 S_s\big((V^HFV)^{-1}\big),\]
i.e., an equivalent form of \eqref{blmeamajorab}
according to Lemma \ref{lm:tangentform}.

\subsection{Proof of Lemma \ref{lm:blmervimajor1}}\label{p:blmervimajor1}

As in the proof of Lemma \ref{lm:blmervi1},
we use the Hermitian positive semidefinite matrix $\widetilde{A}=A-\la_nI_n$
and the Hermitian positive definite matrix
$\widetilde{D}=\diag(\la_1-\la_n,\,\ldots,\,\la_i-\la_n)$
together with their square roots
$\widetilde{A}^{1/2}$ and $\widetilde{D}^{1/2}$.
In addition, we define the orthonormal basis matrix
$\widetilde{Y}=Y'_i\big({Y'_i}^HY'_i\big)^{-1/2}$ of $\Y'_i$
by using the basis matrix $Y'_i=f(A)Y_i$ introduced in Lemma \ref{lm:blmervi}.
Then the Ritz values $\widetilde{\eta}_1\ge\cdots\ge\widetilde{\eta}_i$
of $A$ in $\Y'_i$ coincide with the eigenvalues of $\widetilde{Y}^HA\widetilde{Y}$
so that
\[\Lambda(\widetilde{Y}^H\widetilde{A}\widetilde{Y})
 =\Lambda(\widetilde{Y}^HA\widetilde{Y}-\la_nI_i)
 =[\widetilde{\eta}_1-\la_n,\,\ldots,\,\widetilde{\eta}_i-\la_n]
 \le[\eta'_1-\la_n,\,\ldots,\,\eta'_i-\la_n]\]
holds according to $\Y'_i\subseteq\Y'$ and the Courant-Fischer principles.

Moreover, $X_i^H\widetilde{Y}$ is invertible since
the invertibility of $X_i^HY'_i$ is verified at the beginning of the proof
of Lemma \ref{lm:blmervi}. Then we transform $X_i^H\widetilde{Y}$ as
\[X_i^H\widetilde{Y}=\widetilde{D}^{-1/2}\widetilde{D}^{1/2}X_i^H\widetilde{Y}
 =\widetilde{D}^{-1/2}\big(X_i\widetilde{D}^{1/2}\big)^H\widetilde{Y}
 =\widetilde{D}^{-1/2}\big(\widetilde{A}^{1/2}X_i\big)^H\widetilde{Y}
 =\widetilde{D}^{-1/2}X_i^H\widetilde{A}^{1/2}\widetilde{Y}\]
so that
\[\widetilde{D}^{1/2}=X_i^H\big(\widetilde{A}^{1/2}\widetilde{Y}\big)
 (X_i^H\widetilde{Y})^{-1}.\]
Then applying Lemma \ref{lm:major} to
\[B_1=X_i^H,\quad B_2=\widetilde{A}^{1/2}\widetilde{Y},\quad
 B_3=(X_i^H\widetilde{Y})^{-1},\]
$t=i$ and $c=2$ yields a majorization relation which can be extended as
\[\begin{split}
 &\left[\frac{\la_1-\la_n}{\eta'_1-\la_n},\,\ldots,\,
 \frac{\la_i-\la_n}{\eta'_i-\la_n}\right]\le
 \left[\frac{\la_1-\la_n}{\widetilde{\eta}_1-\la_n},\,\ldots,\,
 \frac{\la_i-\la_n}{\widetilde{\eta}_i-\la_n}\right]
 \\[1ex]
 =\ &\Lambda(\widetilde{D})/\Lambda(\widetilde{Y}^H\widetilde{A}\widetilde{Y})
 =S_i^2(\widetilde{D}^{1/2})/S_i^2(\widetilde{A}^{1/2}\widetilde{Y})
 \prec_w S_i^2(X_i^H)S_i^2\big((X_i^H\widetilde{Y})^{-1}\big)\\[1ex]
 =\ &[1,\ldots,1]\,
 \left[\big(\cos(\theta_1)\big)^{-2},\,\ldots,\,\big(\cos(\theta_i)\big)^{-2}\right]
 =\left[1+\tan^2(\theta_1),\,\ldots,\,1+\tan^2(\theta_i)\right]
\end{split}\]
with the associated principal angles $\theta_1\ge\cdots\ge\theta_i$. Thus
\[\left[\frac{\la_1-\la_n}{\eta'_1-\la_n},\,\ldots,\,
 \frac{\la_i-\la_n}{\eta'_i-\la_n}\right] \prec_w
 \left[1+\tan^2(\theta_1),\,\ldots,\,1+\tan^2(\theta_i)\right]\]
holds and implies bound \eqref{blmerviamajor1}
by subtracting $1$ from each component.
\hfill$\Box$

\begin{remark}\label{mod}
Bound \eqref{blmerviamajor1} is comparable with
\cite[Theorems 2.3 and 2.4]{ka10}
where actually (despite different notations) the smallest Ritz value
of $A$ in $\X_i{\,+\,}\Y'_i$ or $\X_i{\,+\,}\Y'$ is used instead of
the smallest eigenvalue $\la_n$.
Such a modification can be made by restricting
the proof of Lemma \ref{lm:blmervimajor1} to the associated subspace.
In comparison to the corresponding derivation in \cite{ka10},
our new derivation is more direct (e.g., without shifting) and the factorization
$\widetilde{D}^{1/2}=X_i^H\big(\widetilde{A}^{1/2}\widetilde{Y}\big)
 (X_i^H\widetilde{Y})^{-1}$ enables a simpler formulation.
\end{remark}

\def\refname{
\end{document}